\documentclass[a4paper, 11pt]{article}
\usepackage{amsmath, amssymb,amscd, amsthm}
\usepackage[mathscr]{eucal}
\usepackage{graphics}
\usepackage{fullpage}
\newcommand\cyr{%
 \renewcommand\rmdefault{wncyr}%
 \renewcommand\sfdefault{wncyss}%
 \renewcommand\encodingdefault{OT2}%
\normalfont\selectfont} \DeclareTextFontCommand{\textcyr}{\cyr}

\newtheorem{theorem}{Theorem}
\newtheorem{lemma}[theorem]{Lemma}
\newtheorem{corollary}[theorem]{Corollary}

\newtheorem{remark}[theorem]{Remark}

\newtheorem{question}[theorem]{Question}

\newtheorem*{thm5}{Theorem 5}
\newtheorem*{thm6}{Theorem 6}
\newtheorem*{thm9}{Theorem 9}
\newtheorem*{thm7}{Theorem 7}
\newtheorem*{thm10}{Theorem 10}

\def\N{\mathbb N}

\def\T{\mathbb T}

\def\R{\mathbb R}

\begin{document}
\title{Estimates for the Square Variation of Partial Sums of Fourier Series and their Rearrangements}
\author{Allison Lewko\thanks{Supported by a National Defense Science and Engineering Graduate Fellowship.} \; and Mark Lewko}
\date{}
\maketitle
\abstract{We investigate the square variation operator $V^2$ (which majorizes the partial sum maximal operator) on general orthonormal systems (ONS) of size $N$. We prove that the $L^2$ norm of the $V^2$ operator is bounded by $O(\ln(N))$ on any ONS. This result is sharp and refines the classical Rademacher-Menshov theorem. We show that this can be improved to $O(\sqrt{\ln(N)})$ for the trigonometric system, which is also sharp. We show that for any choice of coefficients, this truncation of the trigonometric system can be rearranged so that the $L^2$ norm of the associated $V^2$ operator is $O(\sqrt{\ln\ln(N)})$. We also show that for $p>2$, a bounded ONS of size $N$ can be rearranged so that the $L^2$ norm of the $V^p$ operator is at most $O_p(\ln \ln (N))$ uniformly for all choices of coefficients. This refines Bourgain's work on Garsia's conjecture, which is equivalent to the $V^{\infty}$ case. Several other results on operators of this form are also obtained. The proofs rely on combinatorial and probabilistic methods. }

\section{Introduction}

Let $\T:=[0,1]$ denote the unit interval with Lebesgue measure $dx$ and let $\Phi:=\{\phi_{n}\}_{n \in \N}$ denote an orthonormal system (ONS) of real or complex valued functions on $\T$. By an ONS, we will always mean the set of orthonormal functions $\{\phi_n\}_{n \in \N}$ \emph{and} the ordering inherited from the index set $\N$. For $f \in L^2$, we let $a_{n}= \left<f, \phi_n \right>$ denote the Fourier coefficients of $f$ with respect to the system $\Phi$.  Associated to an ONS is the maximal partial sum operator

\[\mathcal{M}f(x) := \sup_{N}\left|\sum_{n=1 }^{N} a_{n}\phi_{n}(x)\right|.\]

It is well known that the $L^2$ boundedness of the operator $\mathcal{M}$ implies the almost everwhere convergence of the partial sums of the expansion of $f \in L^2$ in terms of the ONS $\Phi$. Almost everywhere convergence is known to fail for some ONS, hence the maximal function $\mathcal{M}$ is known to be an unbounded operator on $L^2$ for some ONS. There is an optimal estimate known for general ONS.

\begin{theorem}\label{RM}\emph{(Rademacher-Menshov)} Let $\{\phi_{n}\}_{n\in\N}=\Phi$ and $f \in L^2$ be as above. Then,

\[ ||\mathcal{M}f||_{L^2} \ll \left(\sum_{n=1}^{\infty}  |a_{n}|^2\ln^2(n+1)\right)^{\frac{1}{2}}\]

where the implied constant is absolute. Moreover, the function $\ln^2(n+1)$ cannot be replaced with any function that is $o(ln^2(n+1))$.
\end{theorem}

This last claim is quite deep and is due solely to Menshov.

While this estimate is optimal in general, it can be improved for many specific systems. For instance, the inequality $||\mathcal{M}f||_{L^2} \ll ||f||_{L^2}$ is known to hold when $\Phi$ is taken to be the trigonometric, Rademacher, or Haar systems. We recall the definitions of these systems in the next section.

Recently, variational norm refinements of the maximal function results stated above have been investigated. To state these results, we first need to introduce some notation. Let $a=\{a_{n}\}_{n=1}^{\infty}$ be a sequence of complex numbers. Then we define the $r$-variation as:
\[||a||_{V^{r}}:= \lim_{K\rightarrow \infty} \sup_{\mathcal{P}_{K}} \left( \sum_{I \in \mathcal{P}_{K}}\left|\sum_{n \in I} a_n \right|^{r} \right)^{1/r},\]
where the supremum is taken over all partitions $\mathcal{P}_K$ of $[K]$ (i.e. all ways of dividing $[K]$ into disjoint subintervals). When $a$ is a finite sequence of length $K$, the quantity is defined by dropping the $\lim_{K\rightarrow \infty}$.

One can easily verify that this is a norm and is nondecreasing as $r$ decreases. Now we will denote the sequence $\{a_n \phi_n(x)\}_{n=1}^{\infty}$ by $S[f](x)$. (Note that this is slightly different than the notation used in \cite{OSTTW}.) When we write $||S[f]||_{V^r}(x)$, we mean the function on $\T$ whose value at $x\in\T$ is obtained by assigning the $r$-th variation of the sequence $S[f](x)$. Furthermore, $||S[f]||_{L^p(V^{r})}$ is the $L^p$ norm of this function. Alternately, we have
\[||f||_{V^{2}}(x) = \sup_{K} \sup_{n_{0}<\ldots<n_{K}}\left(\sum_{l=1}^{K}|S_{n_{l}}[f](x) - S_{n_{l-1}}[f](x)|^{2} \right)^{1/2},\]
where $S_{n_l}[f](x) = \sum_{n=1}^{n_l}a_n\phi_n(x)$ is the $n_l$-th partial sum.

We note that the function $||S[f]||_{V^\infty}(x)$ is essentially the maximal function. More precisely, $\mathcal{M}f(x) \ll ||S[f]||_{V^\infty}(x) \ll \mathcal{M}f(x)$. Since the quantity $||a||_{V^{r}}$ is nondecreasing as $r$ decreases, we see that $||S[f]||_{V^r}(x)$ majorizes the maximal function whenever $r< \infty$.  In \cite{OSTTW}, the following is proved for the trigonometric system $\{e^{2\pi i n x}\}_{n=1}^{\infty}$:

\begin{theorem}\label{varCarleson}Let $r>2$ and $r' < p < \infty$, where $\frac{1}{r} + \frac{1}{r'} = 1$. Then
\[||S[f]||_{L^{p}(V^{r})} \leq C_{p,r} ||f||_{L^p},\]
where $C_{p,r}$ is a constant depending only on $p$ and $r$.
\end{theorem}

This result is rather deep, being a strengthened version of the celebrated work of Carleson and Hunt on the almost everywhere convergence of Fourier series. The analogous inequalities were previously obtained in \cite{JonesWang} in the simpler situation of Ces\`{a}ro partial sums of the trigonometric system. Moreover, the above inequality is known to hold for the Haar system and more generally for martingale differences by Lepingle's inequality, a variational variant of Doob's maximal inequality. In \cite{OSTTW}, it is shown that the condition $r>2$ is necessary in case of the trigonometric system. Our focus here will be to study the case $p=r=2$ for general ONS. In this direction, we prove (closely following the classical proof):

\begin{theorem}\label{main}Let $\Phi$ be an ONS. Then
\begin{equation}\label{var}
||S[f]||_{L^{2}(V^{2})} \ll \left( \sum_{n=1}^{\infty}|a_{n}|^2 \ln^{2}(n+1) \right)^{1/2}.
\end{equation}
If $||\mathcal{M}f||_{L^2} \ll \Delta(N)||f||_{L^2}$ for all $f= \sum_{n=1}^{N}a_n \phi_n$ for some real valued function $\Delta(N)$, then
\begin{equation}\label{Maxrelate}
|||f||_{L^2(V^2)} \ll \left( \sum_{n=1}^{N} \Delta(n)\ln(n+1) |a_n|^2 \right)^{1/2}.
\end{equation}
\end{theorem}

Interestingly, the first inequality strengthens the Rademacher-Menshov theorem stated above, since the right sides are the same (up to implicit constants), yet we have replaced the maximal function with the square variation operator $V^2$ on the left side. Since the $V^2$ operator dominates the maximal operator, this implies the Rademacher-Menshov theorem and the claim that this result is sharp follows from the sharpness of Rademacher-Menshov. This might lead one to think that the two operators behave similarly, however we will see that the $V^2$ operator is much larger than the maximal operator for the classical systems. Theorem \ref{main} can be refined further for certain classes of ONS, see Section \ref{sec7} for discussion of this.

We can apply (\ref{Maxrelate}) to the trigonometric system with $\Delta(N)=O(1)$, the Carleson-Hunt inequality, and obtain the following corollary:

\begin{corollary}\label{varTrig} Let $\{e^{2\pi i n x}\}_{n=1}^{\infty}$ be the trigonometric system. We then have
\begin{equation}\label{trigsys}
||S[f]||_{L^{2}(V^{2})} \ll \left(\sum_{n=1}^{\infty}|a_{n}|^2 \ln(n+1)\right)^{1/2}.
\end{equation}
Moreover, the function $\ln(n+1)$ cannot be replaced by a function that is $o(\ln(n+1))$.
\end{corollary}

The lower bound can be obtained by considering the Dirichlet kernel $D_N(x)=\sum_{n=1}^{N} e^{2 \pi i n x}$. A proof of this is contained in Section 2 of \cite{OSTTW}. Strictly speaking, they work with the de la Vallee-Poussin kernel there, but the same proof works for the Dirichlet kernel.

As we will see below, it is easy to construct an infinite ONS such that $||S[f]||_{L^{2}(V^{2})} \ll ||f||_{L^2}$ holds, by choosing the basis functions $\phi_{n}(x)$ to have disjoint supports. However, this is a very contrived ONS, and it is then natural to ask if there exists a complete ONS such that $||S[f]||_{L^{2}(V^{2})} \ll ||f||_{L^2}$. This is not possible. In fact, we show slightly more:

\begin{theorem}\label{completeDiverg}Let $\{\phi_{n}\}$ be a complete orthogonal system. There exists a $L^{\infty}$ function such that $||S[f]||_{V^{2}}(x)=\infty$ for almost every $x$.
\end{theorem}

In general, this divergence cannot be made quantitative. We show that for any function $w(n)\rightarrow \infty$, there exists a complete ONS such that $||S[f]||_{L^{2}(V^{2})} \ll w(N) ||f||_{L^2}$ whenever $f(x) = \sum_{n=1}^{N}a_n\phi_n(x)$. However, a quantitative refinement is possible if we restrict our attention to uniformly bounded ONS:

\begin{theorem}\label{boundedDiverg} In the case of a uniformly bounded ONS, it is not possible for $w(N)=o(\sqrt{\ln\ln(N)})$. However, there do exist uniformly bounded ONS such that $w(N)=O(\sqrt{\ln\ln(N)})$.
\end{theorem}
The Rademacher system provides an example of the second claim. See Theorem \ref{varRad} below.

Recall that we defined an ONS to be a sequence of orthonormal functions with a specified ordering. This is essential since the behavior of the maximal and variational operators depend heavily on the ordering. For instance, the Carleson-Hunt bound on the maximal function for the trigonometric system makes essential use of the ordering of the system, and the result is known to fail for other orderings. It is thus natural to ask what one can say about the $V^2$ operator for reorderings of the trigonometric system. Surprisingly, it turns out that the $O(\sqrt{\ln (N)})$ bound can be improved to $O(\sqrt{\ln\ln(N)})$ for any choice of coefficients by reordering the system. More generally:

\begin{theorem}\label{mod1Perm}Let $\{\phi_n\}_{n=1}^{N}$ be an ONS such that $|\phi_n(x)|=1$ for all $x$ and $n$, and let $f(x) = \sum_{n=1}^{N} a_n \phi_n(x)$. Then there exists a permutation $\pi:[N]\rightarrow [N]$ such that
\[||f||_{L^2(V^2)} \ll \sqrt{\ln \ln (N)} ||f||_{L^2} \]
holds (for sufficiently large $N$) with respect to the rearranged ONS $\{\psi_n\}_{n=1}^{N}$, where $\psi_n(x) := \phi_{\pi(n)}(x)$.
\end{theorem}

This is perhaps the most technically interesting part of the paper. This result should be compared to Garsia's theorem \cite{Garsia2}, which states that the Fourier series of an arbitrary function with respect to an arbitrary ONS can be rearranged so that the maximal function is bounded on $L^2$. Garsia's proof proceeds by selecting a uniformly random permutation, and arguing that it will satisfy the claim with positive probability. In our case, however, we randomize over a subset of all permutations. This subset is chosen based on structural information about the Fourier coefficients of the function. It is unclear if this restriction is necessary or an artifact of our proof techniques. It would be interesting to extend this result to more general ONS. We note that it can be seen from the work of Qian \cite{Qian} (see also our refinement \cite{Lewko})  that $|| \sum_{n=1}^{N}r_n||_{L^2(V^2)} \gg \sqrt{N\ln\ln(N)}= \sqrt{\ln\ln(N)}\; ||\sum_{n=1}^{N}r_n||_{L^2}$, regardless of the ordering of the Rademacher functions $r_n$, hence the $\sqrt{\ln \ln (N)}$ term in the statement of the theorem is sharp. A similar result can be obtained for general ONS when the coefficients are multiplied by random signs:

\begin{theorem}\label{randomSigns} Let $\{\phi_n\}_{n=1}^{N}$ be an ONS and $f(x) = \sum_{n=1}^{N} a_n \phi_n(x)$. Then there exists a sequence of signs $\epsilon_{n}$ such that
\[||g||_{L^2(V^2)} \ll_{M} \sqrt{\ln \ln (N)} \; ||g||_{L^2} \]
holds, where $g(x) = \sum_{n=1}^{N} \epsilon_n a_n \phi_n(x)$.
\end{theorem}

This easily follows from the following inequality:

\begin{theorem}\label{varRad} Let $\{r_n\}_{n=1}^{N}$ be a sequence of uniformly bounded independent random variables. Then
\[\left|\left| \sum_{n=1}^{N} a_n r_n \right|\right|_{L^2(V^2)} \ll \sqrt{\ln\ln(N)} \left( \sum_{n=1}^{N} a_n^2 \right)^{1/2}.\]
In particular, combining this with Theorem \ref{boundedDiverg}, we see that the $L^2$ norm of the $V^2$ operator for the Rademacher system grows like $\sqrt{\ln\ln(N)}$.
\end{theorem}

Finally, we prove that the $V^p$ norm of some systems can be improved uniformly for all choices of coefficients by a rearrangement, for $p>2$.

\begin{theorem}\label{varVp}Let $\{\phi_{n}\}_{n=1}^{N}$ be an ONS such that $||\phi_{n}||_{L^\infty}\leq M$ for each $n$, and let $p >2$. There exists a permutation $\pi:[N]\rightarrow [N]$ such that the orthonormal system $\{\phi_{\pi(n)}\}_{n=1}^{N}$ satisfies
\begin{equation}\label{varperm}
||S[f]||_{L^{2}(V^{p})} \ll_{M,p} \ln\ln(N)||f||_{L^2}
\end{equation}
for all $f = \sum_{n=1}^{N} a_n \phi_n $.
\end{theorem}

The maximal $V^{\infty}$ version of this result is due to Bourgain \cite{Bour} and represents the best progress known towards Garsia and Kolmogorov's rearrangement conjectures. Our methods rely heavily on those developed in that paper. This also leads us to perhaps the most interesting open problem relating to $V^2$ operators:

\begin{question}Does there exist a permutation $\pi: [N] \rightarrow [N]$ such that the $L^2$ norm of the associated $V^2$ operator on the trigonometric system grows like $o(\sqrt{\ln(N)})$?
\end{question}

Our Theorems \ref{mod1Perm} and \ref{varVp} may be viewed as evidence that this may in fact be possible. It is consistent with our knowledge that one could get growth as slow as $\sqrt{\ln\ln(N)}$. It is known that purely probabilistic techniques in the maximal ($V^\infty$) case can only go as far as Bourgain's bound of $\ln \ln(N)$ (see Remark 2 of \cite{Bour}). Thus, finding a permutation that reduces the growth further (Garsia's conjecture is the assertion that there exists a rearrangement that gets to $O(1)$) would require fundamentally new ideas. However, it is consistent with our current knowledge that the purely probabilistic techniques could get one down to $\ln\ln(N)$ in the $V^2$ case. If true, this will certainly require a much more delicate analysis than the methods used here. Theorem \ref{main} combined with the $V^\infty$ case of the previous theorem does give a bound of $\sqrt{\ln(N)}\ln\ln(N)$ for general bounded ONS for the $V^2$ operator. This is a nontrivial improvement for some systems, but not the most interesting classical systems.

\section{Notation and General Remarks}

We will work with ONS defined on the unit interval $\mathbb{T}$. The underlying space $\mathbb{T}$ plays almost no role in our proofs (the role is similar to that of a probability space in probability theory), and one could replace it with an abstract probability space.

We assume that the ONS is real valued in most of our results. In these cases, one can obtain the same results for complex valued ONS by splitting into real and imaginary parts and applying the arguments to each. The details are routine so we omit them. The proof of Theorem \ref{mod1Perm} is the one place where this requires some care, and thus we work with complex valued functions directly there.

We define the trigonometric system to be the system of complex exponentials $\{e^{2\pi i n x}\}_{n=1}^{\infty}$. Typically the trigonometric system is defined to be the doubly infinite system $\{e^{2\pi i n x}\}_{n=-\infty}^{\infty}$ and the maximal and variational operators are defined with respect to the symmetric partial sums. However, we find it more convenient to define the trigonometric system this way and avoid having to state all of the following results for both singly and doubly infinite systems. All of our results can easily be transferred to the doubly infinite setting (using symmetric partial sums) by splitting the Fourier series of a function $f \in L^2(\T)$ with respect to a doubly infinite system into two functions with singly infinite Fourier series and applying the results in this setting. For instance, note that

\[\mathcal{M}f(x) := \sup_{N}\left|\sum_{n= -N}^{N} a_{n}\phi_{n}(x)\right| \ll \sup_{N}\left|\sum_{n= -N}^{0} a_{n}\phi_{n}(x)\right| + \sup_{N}\left|\sum_{n= 1}^{N} a_{n}\phi_{n}(x)\right|. \]

Thus it follows that the $L^2$ boundedness of the maximal operator associated to the system $\{e^{2\pi i n x}\}_{n=1}^{\infty}$ implies the $L^2$ boudedness of the symmetric maximal operator associated to $\{e^{2\pi i n x}\}_{n=-\infty}^{\infty}$, and similarly for the $V^p$ operators.

The Haar system, which we denote by $\{\mathcal{H}_{n}\}_{n=0}^\infty$, is a complete ONS comprised of the following functions. For $k \in \N$ and $1 \leq j \leq 2^k$, we define $\{\mathcal{H}_{k,j}\}$ by
\[\mathcal{H}_{k,j}(x) = \begin{cases}\sqrt{2^{k}} \quad & x \in \left(\frac{j-1}{2^{k}}, \frac{j-1/2}{2^k}\right),\\
 -\sqrt{2^{k}} & x \in \left(\frac{j-1/2}{2^{k}}, \frac{j}{2^k}\right), \\0&  \mbox{otherwise.}\end{cases}\]
We form the system $\mathcal{H}_{n}$ by ordering the basis functions $\{\mathcal{H}_{k,j}\}$ first by the parameter $k$ and then by the parameter $j$, or $\mathcal{H}_{n}=\mathcal{H}_{j,k}$ for $n=2^{k}+j$. Lastly, we set $H_{0}=1$.

The Rademacher system, denoted $\{r_n(x)\}_{n=1}^{\infty}$, is defined by
\[r_n(x) = \text{sign} \sin \left(2^n \pi x \right).\]
The Rademacher system can also be thought of as independent random variables which take each of the values $\{-1,1\}$ with probability $1/2$.

\section{Variational Rademacher-Menshov-Type Results}\label{sec:rad}
We start by giving a proof of Theorem \ref{main}.

It suffices to assume that $N$ is a power of $2$, say $N=2^{\ell}$. For all $i,k$ such that $0 \leq i \leq \ell$ and $0 \leq k \leq 2^{\ell-i}-1$, we consider the collection of intervals $I_{k,i} := (k 2^{i},(k+1)2^{i}]$.

\begin{lemma}\label{lem:binarydecomp} Any subinterval of $S \subset [0,2^{\ell}]$ can be expressed as the disjoint union of intervals of the form $I_{k,i}$, such as
\begin{equation}\label{Idec}
S = \bigcup_{m} I_{k_{m},i_{m}}
\end{equation}
where at most two of the intervals $I_{k_{m},i_{m}}$ in the union are of each size, and where the union consists of at most $2\ell$ intervals.
\end{lemma}
\begin{proof}Let $S=[a,b]$ and set $i' := \max_{I_{k,i} \subseteq S} i$. It follows that there are at most two intervals of the form $I_{k,i'}$ contained in $S$ (otherwise $S$ would contain an interval of the form $I_{k,i'+1}$). Let $r$ denote the right-most element of the interval with the largest $k$ value satisfying $I_{k,i'} \subseteq S$. Now $b-r$ has a unique binary expansion. It easily follows from this that $(r,b]$ can be written as $[r,b] = \bigcup_{m} I_{k_{m},i_{m}}$ where the union contains only one interval of the form $I_{k_{m},i_{m}}$ of any particular size, and these intervals are disjoint. An analogous argument allows us to obtain a decomposition of this form also for $[a,r']$, where $r'$ is the left-most element of an interval with the smallest $k$ value satisfying $I_{k,i'} \subseteq S$. The lemma follows by taking the union of these two decompositions.
\end{proof}

We now prove

\begin{lemma}\label{varRM1}In the notation above, we have that
\begin{equation}\label{var}
||S[f]||_{L^{2}(V^{2})} \ll  \ln(N)  \left(\sum_{n=1}^{\infty}|a_{n}|^2 \right)^{1/2}.
\end{equation}
\end{lemma}

\begin{proof}
By rounding up to the nearest power of two, we can assume without loss of generality that $N = 2^\ell$ for some positive integer $\ell$ (this change will only affect the constants absorbed by the $\ll$ notation). Now, for each $x$, we have some disjoint intervals $J_1, \ldots, J_b \subseteq [N]$ such that:
\[ ||S[f]||_{V^{2}}(x) = \sqrt{ \sum_{j=1}^b \left(\sum_{n \in J_j} a_n \phi_n (x)\right)^2}.\]
It is important to note that these intervals depend on $x$.

By Lemma \ref{lem:binarydecomp}, each $J_j$ can be decomposed as a disjoint union of the form (\ref{Idec}). In this disjoint union of intervals $I_{k_m, i_m}$, each value of $i_m$ appears at most twice. For each $j$ and $i$, we let $I^j_{i}$ denote the union of the (at most two) intervals in the decomposition of $J_j$ which are of length $2^i$. We then have:
\[ ||S[f]||_{V^{2}}(x) = \sqrt{ \sum_{j=1}^b \left( \sum_{i=0}^\ell \sum_{n \in I^j_i} a_n \phi_n (x)\right)^2}.\]
Applying the triangle inequality for the $\ell^2$ norm, this is:
\[ \leq \sum_{i=0}^\ell \sqrt{\sum_{j=1}^b \left(\sum_{n \in I^j_i} a_n \phi_n(x)\right)^2}.\] Now, since each $I^j_i$ is a union of at most two intervals, this implies:
\begin{equation}\label{pointwise}
    ||S[f]||_{V^{2}}(x) \ll \sum_{i=0}^\ell \sqrt{\sum_{k = 0}^{2^{\ell-i}-1} \left( \sum_{n \in I_{k,i}} a_n \phi_n(x) \right)^2 }.
\end{equation}
Notice that we are now summing over all intervals $I_{k,i}$ for each $i$, regardless of the value of $x$.

We take the $L^2$ norm of both sides of (\ref{pointwise}), and apply the triangle inequality to obtain:
\begin{equation}\label{triangle}
||S[f]||_{L^{2}(V^{2})} \ll \sum_{i=0}^\ell \left|\left| \sqrt{\sum_{k = 0}^{2^{\ell-i}-1} \left( \sum_{n \in I_{k,i}} a_n \phi_n(x) \right)^2 }\right|\right|_{L^2}.
\end{equation}

By linearity of the integral and Parseval's identity, we have that
\[\left|\left| \sqrt{\sum_{k = 0}^{2^{\ell-i}-1} \left( \sum_{n \in I_{k,i}} a_n \phi_n(x) \right)^2 }\right|\right|_{L^2} = \left(\sum_{k=0}^{2^{\ell-i}-1} \sum_{n \in I_{k,i}} a_n^2\right)^{\frac{1}{2}} = \left( \sum_{n=1}^{N} a_n^2\right)^{\frac{1}{2}},\]
for each $i$. Combining this with (\ref{triangle}) and noting that there are $\ll \ln N$ values of $i$, we have:
\[||S[f]||_{L^{2}(V^{2})} \ll \ln(N)  \left(\sum_{n=1}^{\infty}|a_{n}|^2 \right)^{1/2}.\]

\end{proof}

We now define a variant of the function $||S[f]||_{V^2}(x)$ which we will denote by $||S_{\text{L}}[f]||_{V^2}(x)$. For each $x$, we define $S_{\text{L}}[f](x)$ to be the sequence of differences of lacunary partial sums of $f$ at $x$, i.e. $S_{\text{L}}[f](x):=\{S_{2^{0}}[f](x),S_{2^{1}}[f](x)-S_{2^0}[f](x),S_{2^{2}}[f](x)-S_{2^1}[f](x),\ldots \}$. As usual, we let $||S_{\text{L}}[f]||_{V^2}(x)$ denote the 2-variation of this function.

\begin{lemma}\label{longRMvar}In the notation above we have that
\[||S_{\text{L}}[f]||_{L^2(V^2)} \ll \left( \sum_{n=1}^{\infty} \ln(n+1)|a_{n}|^2 \right)^{1/2}.\]
\end{lemma}
\begin{proof}
We will need the inequality $|a|^2 \leq 2|a-b|^2 + 2|b|^2$ for any real numbers $a,b$. For each $x$, there exists some sequence $m_0(x), m_1(x), m_2(x), \ldots$ such that:
\begin{equation}\label{sequence}
    ||S_{\text{L}}[f]||^2_{V^2}(x) = \left| S_{2^{m_0(x)}} [f](x)\right|^2 + \sum_{i=1}^{\infty} \left| S_{2^{m_i(x)}}[f](x) - S_{2^{m_{i-1}(x)}}[f](x)\right|^2.
\end{equation}

Setting $a := S_{2^{m_i(x)}}[f](x) - S_{2^{m_{i-1}(x)}}[f](x)$ and $b:= f(x)- S_{2^{m_{i-1}(x)}}[f](x)$, we can apply the inequality above to obtain:
\[\left| S_{2^{m_i(x)}}[f](x) - S_{2^{m_{i-1}(x)}}[f](x)\right|^2 \leq 2\left| S_{2^{m_i(x)}} [f](x) - f(x)\right|^2 + 2\left|S_{2^{m_{i-1}(x)}}[f](x) - f(x)\right|^2\]
for each $i\geq 1$. Combining this with (\ref{sequence}), we have:
\begin{eqnarray*}
  \nonumber ||S_{\text{L}}[f]||^2_{V^2}(x) &\ll & \left| S_{2^{m_0(x)}} [f](x)\right|^2  + \sum_{i=1}^\infty \left| S_{2^{m_i(x)}} [f](x) - f(x)\right|^2 + \left|S_{2^{m_{i-1}(x)}}[f](x) - f(x)\right|^2\\
  \nonumber &\ll & \left| S_{2^{m_0(x)}} [f](x)\right|^2 + \sum_{i=0}^\infty \left| S_{2^{m_i(x)}} [f](x) - f(x)\right|^2 \\
  \nonumber & \ll & \left| S_{2^{m_0(x)}} [f](x)\right|^2 + \sum_{m=0}^\infty \left| S_{2^m} [f](x) - f(x) \right|^2.
\end{eqnarray*}
Note that in this last quantity, we are always summing over all values of $m$, instead of summing over a subsequence dependent on $x$.

This gives us
\[||S_{\text{L}}[f]||_{V^2}(x) \ll \left(\left| S_{2^{m_0(x)}} [f](x)\right|^2 + \sum_{m=0}^\infty \left| S_{2^m} [f](x) - f(x) \right|^2\right)^{\frac{1}{2}}.\]
Now we take the $L^2$ norm of both sides of this inequality to obtain:
\[||S_{\text{L}}[f]||_{L^2(V^2)} \ll \left( \sum_{n=1}^\infty \ln(n+1) a_n^2\right)^{\frac{1}{2}}.\]
To see this, note that $\left|S_{2^m} [f](x) -f(x) \right| = \left| \sum_{n=2^m+1}^\infty a_n \phi_n(x)\right|$ and each $n$ is greater than $2^m$ for $\ll \ln(n)$ values of $m$. The result then follows from Parseval's identity.

\end{proof}

We now combine these two results to prove the following theorem.

\begin{theorem}\label{StrongRMvar}For an arbitrary ONS, in the notation above, we have
\[||S[f]||_{L^2(V^2)} \ll \left( \sum_{n=1}^\infty \ln^2 (n+1) a_n^2\right)^{\frac{1}{2}}.\]
\end{theorem}

\begin{proof}
We write $U_{k}(x):= \sum_{n= 2^{k-1}+1}^{2^{k}} a_n \phi_n(x)$ (when $k=0$, $U_0 (x) := a_1 \phi_1(x)$.). We claim that
\[ ||S[f]||_{L^2(V^2)}^2 \ll \int_{\T} \left( ||S_{\text{L}}[f]||_{V^2}^2(x) + \sum_{k=0}^{\infty} || U_k||_{V^2}^2(x) \right) dx. \]

To see this, note that any interval $[a,b]$ can be decomposed as the disjoint union of at most three intervals $I_l, I_c, I_r$, where $I_c=(2^{k},2^{k'}]$ and $I_l \subseteq (2^{k-1},2^{k}]$ and $I_r \subseteq (2^{k'},2^{k'+1})$ (here, $2^k$ can be set as the smallest integral power of 2 contained in $[a,b]$, and $2^{k'}$ can be set as the largest integral power of 2 contained in $[a,b]$). Now,
 $\int_{\T} ||S_{\text{L}}[f]||_{V^2}^2(x)  dx \ll \sum_{n=1}^{\infty} \ln(n+1)|a_{n}|^2 $ from the previous lemma, which is clearly bounded by $ \sum_{n=1}^\infty \ln^2 (n+1) a_n^2$. By Lemma \ref{varRM1}, we have
\[\int_{\T}  || U_k||_{V^2}^2(x) dx \ll \ln^2(2^{k}+1) \sum_{n= 2^{k-1}+1}^{2^k}  a_n^2  \ll \sum_{n= 2^{k-1}+1}^{2^k} \ln^2(n+1)  a_n^2.     \]
Combining these estimates completes the proof.
\end{proof}

Next we show that these estimates can be improved if one has additional information regarding the ONS. In particular, if the partial sum maximal operator $\mathcal{M}$ associated to the system is bounded then one can replace the $\ln^2(n)$
above with an $\ln(n)$.

\begin{theorem}Let $f(x) = \sum_{n=1}^{N}a_n \phi_n(x)$ and assume that $||\mathcal{M}f||_{L^2}  \ll \Delta(N) \left( \sum_{n=1}^{N}a_n^2 \right)^{1/2}$ for any choice of $f$. Then
\[ ||f||_{L^2(V^2)} \ll \Delta(N) \sqrt{\ln(N)} \left( \sum_{n=1}^{N}a_n^2 \right)^{1/2} \]
and
\[ ||f||_{L^2(V^2)} \ll \left( \sum_{n=1}^{N} \Delta(n)\ln(n+1) a_n^2 \right)^{1/2}.\]
In particular, if the quantity on the right is finite, then the variational operator applied to $f$ must be finite almost everywhere.
\end{theorem}

\begin{proof} As before, without loss of generality, we may assume that $N=2^\ell$ for some positive integer $\ell$. And we consider the collection of dyadic subintervals of $[1,N]$ of the form $I_{k,i} = (k 2^{i},(k+1)2^{i}]$ for each $0 \leq i \leq \ell$, $0 \leq k \leq 2^{\ell-i}-1$.
We will refer to intervals of this form as admissible intervals.

Now we note that an arbitrary interval $J=[a,b] \subseteq [N]$ can be written as a disjoint union $J = J_l \cup J_r$, where $J_r \subseteq I_{{k_r},{i_r}}$ and $J_l \subseteq I_{{k_l},{i_l}}$ and $|J_l| \geq \frac{1}{2}|I_{{k_l},{i_l}}|$ and
$|J_r| \geq \frac{1}{2}|I_{{k_r},{i_r}}|$. We allow one of the intervals to be empty if needed, although in the following we will always assume that the intervals are not empty, since estimating the contribution from an empty interval is trivial.
That is, we can write an arbitrary interval $J$ as the union of two intervals which are contained within admissible intervals and the intersection with the admissible intervals is a constant fraction of the the admissible interval.

For $J \subseteq [N]$, let $S_J := \sum_{n \in J} a_n \phi_n(x)$. We now claim the pointwise inequality
\[ ||f||_{V^2}^{2}(x) \ll  \sum_{0 \leq i \leq \ell} \; \sum_{0 \leq k \leq 2^{\ell-i}-1} |\mathcal{M}S_{I_{k,i}}(x)|^2.\]

Note that the sum on the right is only over all admissible intervals. To see that this inequality holds, let $ \{J_{i}\}_{i=1}^{m}$ be a partition of $[N]$ that maximizes the square variation (at $x$). From the discussion above,
we can associate disjoint $J_i^{l}$ and $J_i^{r}$ to $J_i$ such that $J_i \subset J_i^{l} \cup J_i^{r}$. Moreover, we can find disjoint admissible intervals $I_i^{l}$ and $I_i^{r}$ such that $J_{i}^{s} \subseteq  I_i^{s}$ and
$|J_i^s| \geq \frac{1}{2} |I_i^{s}|$ ($s \in \{r,l\}$).

We observe that $|S_{J_i}(x)|^2 \ll |\mathcal{M}S_{I_i^{l}}(x)|^2 + |\mathcal{M}S_{I_i^{r}}(x)|^2$. Moreover, any particular admissible interval $I$ will be associated to at most two intervals in the partition $\{J_i\}$ since
the intervals in the partition are disjoint and have at least half the length of the associated admissible interval. The pointwise inequality above now follows. Now integrating each side, applying the hypothesized inequality
$||\mathcal{M}S_J||_{L^2}^2  \ll \Delta^2(N) \sum_{n \in J}a_n^2$, and noting that every point in $[N]$ is in $O(\ln(N))$ admissible intervals, we have that
\[\int_{\T}||f||_{V^2}^{2} dx \ll\sum_{0 \leq i \leq \ell}\;  \sum_{0 \leq k \leq 2^{\ell-i}-1} \int_{\T}|\mathcal{M}S_{I_{k,i}}(x)|^2 dx \]
\[ \ll \Delta^2(N) \ln(N) \sum_{n=1}^{N} a_n^2.\]

Taking the square root of each side completes the the proof of the first inequality in the theorem statement. The second statement follows from the first via the argument used to prove Theorem \ref{StrongRMvar}. Note that we obtained
a bound on the lacunary partial sums in Lemma \ref{longRMvar} of the order $\sqrt{\ln(n)}$. This estimate was better than we needed for the proof of Theorem \ref{StrongRMvar}, however is exactly the order we need here.
\end{proof}

This completes the proof of Theorem \ref{main} and Corollary \ref{varTrig} follows.

\section{Lower bounds}

In this section, we prove:

\begin{thm5} Let $\{\phi_n(x)\}$ be a complete ONS. Then there exists a function $f \in L^{\infty}(\mathbb{T})$ such that for almost every $x \in \T$
\begin{equation}\label{fail}
||f||_{V^2}(x) = \infty.
\end{equation}
\end{thm5}

Here, as before, $||f||_{V^{2}}(x) = \sup_{K} \sup_{n_{0}<\ldots<n_{K}}\left(\sum_{l=1}^{K}|S_{n_{l}}[f](x) - S_{n_{l-1}}[f](x)|^{2} \right)^{1/2}$ where $S_{n_l}[f](x) = \sum_{n=1}^{n_l}a_n\phi_n(x)$ is the $n_l$-th partial sum.

Using Lemma \ref{rw} below and properties of the Dirichlet kernel, Jones and Wang showed (\ref{fail}) for the trigonometric system. In the case of general orthonormal systems, we do not have analytic information regarding the partial summation operator and need to proceed differently.  We start by establishing the result for the Haar system.

We let $E_{k}:L^1 \rightarrow L^1$ denote the conditional expectation operator defined as follows. For $x \in [ l 2^{-k}, (l+1) 2^{-k} )$, $0 \leq l < 2^{k}$, $l \in \N$ we define \[E_{k}f(x) = \int_{l 2^{-k}}^{(l+1) 2^{-k}}f(x)dx.\]

Using a probabilistic result of Qian \cite{Qian}, Jones and Wang \cite{JonesWang} showed that:

\begin{lemma}\label{rw} (Proposition 8.1 of \cite{JonesWang}) There exists $f\in L^{\infty}(\T)$ such that
\[\sup_{K} \sup_{n_{0}<\ldots<n_{K}}\left(\sum_{\ell=1}^{K}|E_{n_{\ell}}f(x) - E_{n_{\ell-1}}f(x)|^{2} \right)^{1/2} = \infty\]
 almost everywhere.
\end{lemma}

If we let $S_n[f]$ denote the partial summation operator with respect to the Haar system, then it easily follows that $E_kf(x) =S_{n_{k+1}}[f](x)- S_{n_{k}}[f](x)$ for some sequence $\{n_k\}$.  Therefore, there exists $f \in L^{\infty}(\T)$ such that $||f||_{V^2}(x)= \infty$ for almost every $x \in \T$, where the operator $V^2$ is associated to the Haar system. For future use, let us define $\{b_{n}\}$ to be the Haar coefficients of the function $f$, that is
\begin{equation}\label{bdef}
b_{n} = \left<f(x), \mathcal{H}_{n}(x) \right>.
\end{equation}

We will also need a theorem of Olevskii (see \cite{Olev} Chapter 3), which requires that we introduce some additional notation. Let $\{g_n\}$ and $\{f_n\}$ be two sequences of real-valued measurable functions on $\T$. We say that they are weakly isomorphic if for each $n \in \N$ there exists an invertible measure-preserving mapping $T_{n}:\T \rightarrow \T$ that is one-to-one on a set of full measure and satisfies

\[f_{k}(T_n x ) = g_{k}(x)\]

for all $1\leq k \leq n$.

\begin{theorem}\label{thm:olevskii} (Olevskii) Let $\{\phi_{n}\}_{n=1}^{\infty}$ be a complete real-valued orthonormal system. There exists an orthonormal system $\{H_{k}\}_{k=1}^{\infty}$ that is weakly isomorphic to the Haar system, and a sequence $\{n_{k}\}_{k=1}^{\infty}$ such that
\[ \left|\left|\sum_{i=n_{k}+1}^{n_{k+1}}  \left< H_{j}, \phi_{i} \right> \phi_{i}(x)\right|\right|_{L^2} \leq 2^{-k-j}\]
whenever $j \neq k$.
\end{theorem}

We now set $ \tilde{f}(x) := \sum_{n=1}^{\infty}b_n H_{n}(x)$, for $b_n$ defined in (\ref{bdef}). Using the fact that the (finite) partial sums of the series defining $\tilde{f}(x)$ are weakly isomorphic to the partial sums of the Haar expansion of $f$, it follows that the partial sums of the function $\tilde{f}$ are uniformly bounded, hence $\tilde{f} \in L^{\infty}(\T)$.

\begin{lemma}For $\tilde{f}$ defined as above, we set $c_{n}:= \left<\tilde{f},\phi_n\right>$. It follows that
\[\sum_{n= n_{k}+1}^{n_{k+1}} c_{n}\phi_n(x) = b_{k} H_{k}(x) + e_{k}(x), \]
where $\sum_{k}|e_{k}(x)| < \infty$ for almost every $x$.
\end{lemma}

\begin{proof}
Since $\tilde{f}(x)= \sum_{j=1}^{\infty} b_j H_j (x)$, we have
\[\sum_{n_{k}+1}^{n_{k+1}} c_{n}\phi_n(x) = \sum_{n=n_{k}+1}^{n_{k+1}} \left<\sum_{j=1}^{\infty} b_j H_j (x) , \phi_n(x) \right>\phi_n(x) \]
\[= \sum_{n=n_{k}+1}^{n_{k+1}} b_k\left< H_k (x) , \phi_n(x) \right>\phi_n(x) + \sum_{n=n_{k}+1}^{n_{k+1}} \left<\sum_{j \neq k} b_j H_j (x) , \phi_n(x) \right>\phi_n(x).  \]
By applying the triangle inequality, we obtain:
\[ \left|\left|b_k H_k(x) - \sum_{n_{k}+1}^{n_{k+1}} c_{n}\phi_n(x) \right|\right|_{L^2} \leq |b_k| \left|\left| \sum_{n \notin [n_k+1,n_{k+1}]} \left< H_k (x) , \phi_n(x) \right>\phi_n(x)  \right| \right|_{L^2} \]
\[+ \sum_{j \neq k}|b_j| \left| \left| \sum_{n=n_{k}+1}^{n_{k+1}} \left< H_j (x) , \phi_n(x) \right>\phi_n(x) \right| \right|_{L^2}. \]

Now applying Theorem \ref{thm:olevskii}, we have that
\[\left|\left|b_k H_k(x) - \sum_{n_{k}+1}^{n_{k+1}} c_{n}\phi_n(x) \right|\right|_{L^2} \ll 2^{-k} \left( |b_k| \sum_{j\neq k} 2^{-j} +\sum_{j\neq k} |b_j| 2^{-j} \right) \ll 2^{-k}||\tilde{f}||_{L^2}. \]
The last bound follows from the fact that $|b_j| \leq ||\tilde{f}||_{L^2} = \left(\sum_{i=1}^{\infty} b_i^2\right)^{1/2}$ for all $j$.

Denoting the expression on the inside of the norm on the left as $e_k(x)$, we see that $\big|\big|\sum_{k=1}^{\infty} |e_k| \big|\big|_{L^2} \ll ||\tilde{f}||_{L^2}$ and hence $\sum_{k=1}^{\infty}|e_k(x)|$ is finite for almost every $x \in \mathbb{T}$.

\end{proof}

We now prove Theorem \ref{completeDiverg}. We let $V_{\phi}$ and $V_{H}$ denote the variation operators associated to the systems $\{\phi_{n}\}$ and $\{H_{n}\}$ respectively. Moreover, we let $V^2$ be the variation operator associated to the partial sums of the absolutely convergent function $E(x)=\sum_{k=1}^{\infty}e_{k}(x)$. We have, for almost every $x \in \T$,
\[|| E||_{V^2}(x) \leq \sum_{k=1}^{\infty} |e_{k}(x)| < \infty.\]
It follows that
\[ ||\tilde{f}||_{L^2(V_H^2)} =  \left|\left| \sum_{k=1}^{\infty} b_k H_{k} \right|\right|_{L^2(V_H^2)} \leq ||\tilde{f}||_{L^2(V_\phi^{2})} - ||E||_{L^2(V^2)}.  \]
Since the first quantity in this expression is infinite almost everywhere, and the third quantity is finite almost everywhere, it must hold that $||\tilde{f}||_{L^2(V_\phi^2)}$ is infinite almost everywhere.
This completes the proof of the theorem.

Our proof of Theorem \ref{completeDiverg} was purely qualitative, a feature we inherit from Theorem \ref{thm:olevskii}, which relies on the Riemann-Lebesgue lemma. Next we show that it is impossible to obtain a quantitative lower bound on the growth of the variation in Theorem \ref{completeDiverg}.

\begin{remark} One could obtain the conclusion of Theorem \ref{completeDiverg} for functions in more restrictive classes. Combining the above argument with known perturbation techniques, one can show that the $f$ in the statement of the theorem can be taken to be continuous. The proof of this relies on the fact that one already has an example in $L^{\infty}$ (an example in $L^2$ is not sufficient). See \cite{Olev} p.67 and the associated references for details. Additionally, one can show that for any nonconstant function $f$, there exists an invertible measure preserving transformation of $T:\T \rightarrow \T$ such that the conclusion holds for $g(x) = f(T(x))$. See \cite{Olev} p.69 and the related references for details. From this, we see that one cannot hope to prove that $V^2$ is bounded on $L^2$ even in ``restricted weak type" form, at least not for complete systems. Since the details of these arguments are not essential to our current investigation, and are essentially a combination of the above argument and the ideas of the cited papers, we omit them.
\end{remark}

\begin{theorem}\label{V2bex}Let $w(\cdot)$ denote a positive real-valued function monotonically increasing to infinity. Then there exists a complete orthonormal system $\{\phi_n\}_{n=1}^{\infty}$ such that for all sufficiently large $N \in \N$,
\[ || f ||_{L^2(V^2)} \ll  w(N) \left(\sum_{n=1}^{N} |a_n|^2 \right)^{\frac{1}{2}} .\]
for all $f$ of the form $f(x) = \sum_{n=1}^{N} a_n \phi_n(x)$.

\end{theorem}

\begin{proof} Our example will be a rearrangement of the Haar system. We let $\Psi =\{\psi_n(x)\}_{n=1}^{\infty}$ be a subsequence of the Haar system with disjoint supports. We let $\{\rho_{n}(x)\}_{n=1}^{\infty}$ denote the subsequence of the Haar system consisting of all the elements of the Haar system that are not included in $\Psi$.  We now form a complete orthonormal system $\{\phi_{n}\}$ by sparsely inserting elements of the sequence $\{\rho_{n}(x)\}_{n=1}^{\infty}$ into the sequence $\{\psi_n(x)\}_{n=1}^{\infty}$, maintaining the relative ordering of each sequence. Clearly we may do this so that the first $N$ elements of the system $\{\phi_n\}$ have at most $w(n)$ elements from the $\rho$'s. We thus may partition the indices $[N]$ of the system $\{\phi_{n}\}_{n=1}^{N}$ into two classes. We let $S$ be the subset of indices $n$ for which $\phi_n = \rho_m$ for some $m$ and $S^c := [N] \setminus S$. We note that for $n \in S^c$, $\phi_n$ is an element of the subsequence $\Psi$, and so all of these have disjoint supports.

We then have:
\[\left|\left| \sum_{n \in S} a_n \phi_n + \sum_{n \in S^{c}} a_n \phi_n \right|\right|_{L^2(V^2)} \leq \left|\left| \sum_{n \in S} a_n \phi_n\right|\right|_{L^2(V^2)} + \left|\left|\sum_{m \in S^{c}} a_m \phi_m \right|\right|_{L^2(V^2)} \]
\[\ll \ln(w(n))||f||_{L^2} + ||f||_{L^2}  \ll \ln(w(n))||f||_{L^2} \ll  w(n)||f||_{L^2}.\]
Here, we have employed the triangle inequality, Lemma \ref{varRM1}, and the fact that $\{\phi_n\}_{n \in S^c}$ have disjoint supports.

\end{proof}

Lastly, we show that if a system is uniformly bounded, then an quantitative lower bound on the growth of the $V^2$ operator is available, even without assuming completeness.

\begin{thm6} Let $\{\phi_{n}\}_{n=1}^N$ be an ONS uniformly bounded by $M$. Then there exists a function of the form $f=\sum_{n=1}^{N} a_n \phi_n (x)$ such that
\[ ||S[f]||_{L^2(V^2)} \gg_{M} \sqrt{\ln\ln(N)}||f||_{L^2} \]
\end{thm6}
In light of Theorem \ref{varRad}, this is best possible.

To prove this, we will rely on the following lemma:
\begin{lemma}\label{lem:Qian} We let $c_1, \ldots, c_N$ denote real numbers, all $\geq \delta$ for some constant $\delta >0$. We let $X_1, \ldots, X_N$ denote independent Gaussian random variables, each with mean 0 and variance 1. Then
\[\mathbb{E}\left[ \left|\left|\sum_{n=1}^N c_n X_n \right|\right|_{V^2}\right] \gg \delta \sqrt{N \ln \ln (N)}.\]
\end{lemma}

\begin{proof} We essentially follow the proof of Theorem 2.1 in \cite{Qian} (pp. 1373-1375), with minor modifications. We let $\Phi(x)$ denote the standard normal distribution function. By Lemma 2.1 of \cite{Qian} (p. 1373), we have that
\begin{equation}\label{Qian}
1-\Phi(x) \geq (1/12) exp(-3x^2/4) \text{ for } x\geq 1.
\end{equation}
We define $S_k = \sum_{n=1}^k c_n X_n$ and we set $K:= 25$. We also set
\[\ell := \ell(N):= \left\lfloor \frac{\ln N}{4 \ln K}\right\rfloor \text{ and } m:= m(N) := \left\lfloor \frac{\ln N}{2 \ln K}\right\rfloor.\]
We let $L x:= \max\{1, \ln x\}$.

For each $\omega \in \Omega$ (where $\Omega$ denotes the probability space), we define $E_N(\omega)$ to be the subset of values $t \in \{1, 2, \ldots, N-\sqrt{N}\}$ such that, for some $\ell \leq j \leq m$, $|S_{t+K^j}(\omega) - S_t(\omega)| \geq \delta \sqrt{K^j LL(N)}/2$. Additionally, for each fixed $t$ and $j$, we define the event
\[E_N^j(t):= \left\{\omega: |S_{t+K^j}(\omega)-S_{t+K^{j-1}}(\omega)|\geq \delta \sqrt{K^j LL(N)}\right\}.\]
Now, $S_{t+K^j} -S_{t+K^{j-1}}$ is distributed as a Gaussian random variable with mean 0 and variance equal to
\[\sigma^2 := Var[S_{t+K^j} -S_{t+K^{j-1}}] = \sum_{n=t+K^{j-1}+1}^{t+K^j} c_n^2.\]
For any $\lambda \in \mathbb{R}$,
\[\mathbb{P} \left[ S_{t+K^j}(\omega)-S_{t+K^{j-1}} \geq \lambda\right] = 1 - \Phi \left(\frac{\lambda}{\sigma}\right).\]
We apply this with $\lambda := \delta \sqrt{K^j LL(N)}$, and since each $c_n \geq \delta$, we have:
\[\frac{\lambda}{\sigma} \leq \sqrt{\frac{K^j LL(N)}{K^j - K^{j-1}}}.\]
Therefore, using (\ref{Qian}), we obtain:
\[\mathbb{P}[E_N^j(t)] = 1 - \Phi \left(\frac{\lambda}{\sigma}\right) \geq 1 - \Phi \left(\sqrt{\frac{K^j LL(N)}{K^j - K^{j-1}}}\right) \geq \frac{1}{12} exp\left( -\frac{3}{4} \frac{K^j}{K^j - K^{j-1}} LL (N)\right).\]
This is $\geq \frac{1}{12} exp \left(-\frac{4}{5} LL (N)\right) = \frac{1}{12} (\ln (N))^{-4/5}$.

We observe that if $|S_{t+K^j}(\omega) - S_{t+ K^{j-1}}(\omega)| \geq \delta \sqrt{K^j LL(N)}$ for some $\ell < j \leq m$, then either $|S_{t+K^j}(\omega) - S_t (\omega)| \geq \delta \sqrt{K^j LL(N)}/2$ or $|S_{t+K^{j-1}}-S_t| \geq \delta \sqrt{K^j LL(N)}/2 \geq \delta \sqrt{K^{j-1}LL(N)}/2$. Thus,
\[\omega \in \bigcup_{j=\ell+1}^m E_N^j(t) \Rightarrow  t \in E_N(\omega).\]
Therefore, for any $t \in \{1, 2, \ldots, N -\lfloor \sqrt{N}\rfloor\}$, we have:
\[\mathbb{P}\left[ \omega: t \in E_N(\omega)\right] \geq \mathbb{P}\left[ \bigcup_{j = \ell+1}^m E_N^j(t)\right].\]

We note that for $j' \neq j$, $E_N^j (t)$ and $E_N^{j'}(t)$ depend on disjoint sets of the random variables $X_i$, and so are independent events. Therefore, letting $\overline{E}_N^j(t)$ denote the complement of $E_{N}^j(t)$, we have
\[\mathbb{P}\left[ \bigcup_{j = \ell+1}^m E_N^j(t)\right] = 1 - \mathbb{P}\left[ \bigcap_{j = \ell+1}^m \overline{E}_N^j(t)\right] = 1 - \prod_{j=\ell+1}^m \mathbb{P}[\overline{E}_N^j(t)].\]
By the above computations, this is
\[\geq 1 - exp\left( -(1/12)(m-\ell)(\ln N)^{-4/5}\right).\]
For sufficiently large $N$, we can bound this by:
\[> 1 - exp \left( - (\ln N)^{1/5}/(52 \ln K)\right):= 1 - p_N.\]

This shows that for each $t$, $\mathbb{P}\left[ \omega: t \in E_N(\omega)\right] > 1 - p_N$. We can alternately express this as:
\[\int_{\Omega} 1_{E_N}(t) d\mathbb{P} > 1- p_N,\]
where $1_{E_N}(t)$ denotes the function that is equal to 1 when $t \in E_N(\omega)$ and equal to 0 otherwise. We define the subset $\mathcal{S} \subseteq \Omega$ to be the set of $\omega \in \Omega$ such that $|E_N(\omega)| > (1-\sqrt{p_N})(N-\sqrt{N})$. Then
\begin{equation}\label{Qian2}
\mathbb{P}[\mathcal{S}] > 1 - \sqrt{p_N}.
\end{equation}
To see this, observe that
\[\int_{\Omega} \sum_{t=1}^{N-\sqrt{N}} 1_{E_N}(t) d\mathbb{P} = \sum_{t=1}^{N-\sqrt{N}} \int_{\Omega}1_{E_N}(t) d \mathbb{P} > (N - \sqrt{N})(1 - p_N).\]
Now, if $\mathbb{P}[\mathcal{S}] \leq 1 - \sqrt{p_N}$ held, this would imply that the integral on the left hand side of the above is also
\[\leq \sqrt{p_N}\left(1- \sqrt{p_N}\right)\left(N - \sqrt{N}\right) + \left(1-\sqrt{p_N}\right)\left(N - \sqrt{N}\right) = \left(N - \sqrt{N}\right)\left(1-p_N\right),\]
which is a contradiction.

We next use the following Vitali covering lemma:
\begin{lemma}\label{Vitali} (\cite{Folland}, Lemma 3.15) Let $\mu(A)$ denote the Lebesgue measure of a set $A \subseteq \mathbb{R}$. Let $\mathcal{U}$ be a collection of open intervals in $\R$ with bounded union $W$. Then for any $\lambda < \mu(W)$, there is a finite, disjoint subcollection $\{V_1, V_2, \ldots, V_q\} \subseteq \mathcal{U}$ such that $\sum_{i=1}^q \mu(V_i) \geq \lambda/3$.
\end{lemma}

For sufficiently large $N$, (\ref{Qian2}) implies that with probability $> 1- \sqrt{p_N}$, for $\geq N':= \lfloor (1-\sqrt{p_N})(N-\sqrt{N}-1)\rfloor$ integers $t \in \{1,2, \ldots, N-\sqrt{N}\}$ (we will call them $t_1, t_2, \ldots, t_{N'}$), we have corresponding values $j_1, \ldots, j_{N'}$ (all $\leq m$) such that $|S_{t_i + K^{j_i}} - S_{t_i}| \geq \delta \sqrt{K^{j_i} LL (N)}/2$ for each $i$ from 1 to $N'$. We consider the collection $\mathcal{U}$ of the open intervals $(t_i, t_i + K^{j_i})$ for $i$ from 1 to $N'$. We note that each $K^{j_i} > 1$. We fix some positive constant $\alpha < 1$. For $N$ sufficiently large, we have $N' > \alpha N$. (Note that $p_N$ approaches 0 as $N$ goes to infinity). Therefore, the union of the intervals in $\mathcal{U}$ is a subset of $(0,N]$ with Lebesgue measure $\geq N' > \alpha N$.

Applying Lemma \ref{Vitali}, we conclude that there is disjoint subcollection of these open intervals, denoted by $\{(t_i, t_i + K^{j_i})\}_{i \in Q}$, where $Q \subseteq [N']$, such that
\[\sum_{i \in Q} K^{j_i} \geq \alpha N/3.\]
The closures of the intervals in $Q$ are non-overlapping except for possibly at their endpoints. Relabeling the $t_i$'s for $i \in Q$ as $t_1, \ldots, t_q$ (where $q = |Q|$), we have $t_1 < t_1 + K^{j_1} \leq t_2 < t_2 + K^{j_2} \leq \cdots \leq t_q < t_q + K^{j_q} \leq N$. Then,
\[\sum_{i=1}^q \left(S_{t_i+K^{j_i}} - S_{t_i}\right)^2 \geq (1/4) \delta^2 \sum_{i=1}^q K^{j_i} LL(N) \geq (\alpha/12)\delta^2 N LL(N).\]

This implies that
\[\mathbb{P}\left[ \left|\left| \sum_{n=1}^n c_nX_n\right|\right|_{V^2} \geq \delta \sqrt{(\alpha/12) N \ln \ln N}\right] > 1 - \sqrt{p_N},\]
for all sufficiently large $N$.
Hence, by Markov's inequality,
\[ \mathbb{E}\left[ \left|\left|\sum_{n=1}^N c_n X_n \right|\right|_{V^2}\right] \geq \delta \sqrt{(\alpha/12) N \ln \ln N} (1- \sqrt{p_N}) \gg \delta \sqrt{N \ln \ln N}.\]

\end{proof}

We now prove Theorem \ref{boundedDiverg}. We begin by noting that for each $n$, $\int_{\T} \phi_n^2(x) dx = 1$ and $|\phi_n(x)| \leq M \; \forall x$ implies that there are positive constants $\epsilon, \delta >0$ (depending on $M$) such that for some sets $U_n \subseteq \T$ each of measure $\geq \epsilon$, $|\phi_n(x)| \geq \delta$ for all $x \in U_n$. For each $n$, we let $\chi_n$ denote the characteristic function of the set $U_n$. We then have:
\begin{equation}\label{characteristic}
\int_{\T} \sum_{n=1}^{N} \chi_n(x) dx = \sum_{n=1}^{N} \int_{T} \chi_n(x) dx \geq N\epsilon.
\end{equation}
We define $\epsilon' :=\frac{\epsilon}{2}$. Then the function $\sum_{n=1}^N \chi_n (x)$ must be $\geq \epsilon' N$ on a set of measure $\geq \epsilon'$. To see this, note that $0 \leq \sum_{n=1}^N \chi_n(x)\leq N$ for all $N$. If this function is less than $\epsilon' N$ on a set of measure $> 1 - \epsilon'$, this would imply
\[\int_T \sum_{n=1}^N \chi_n(x) dx < \epsilon' N (1-\epsilon') + \epsilon' N = (1-\epsilon/4) N\epsilon,\]
contradicting (\ref{characteristic}). Thus, there is some set $U$ of measure $\geq \epsilon'$ such that for every $x \in U$, $|\phi_n(x)| \geq \delta$ for at least $\epsilon' N$ values of $n$.

We let $X_1, \ldots, X_N$ denote independent Gaussian random variables with mean 0 and variance 1. We consider the quantity
\[\mathbb{E}\left[ \left|\left| \{X_n \phi_n(x)\}_{n=1}^N \right|\right|^2_{L^2(V^2)}\right].\]
This can be written as:
\[ \mathbb{E}\left[ \int_{\T} \left|\left| \{X_n\phi_n(x)\}_{n=1}^N \right|\right|^2_{V^2} dx \right] = \int_{\Omega} \int_{\T} \left|\left| \{X_n\phi_n(x)\}_{n=1}^N \right|\right|^2_{V^2} dx d\mathbb{P}.\]
By Fubini's theorem, we may exchange the integrals to obtain
\[ = \int_{T} \int_{\Omega} \left|\left| \{X_n\phi_n(x)\}_{n=1}^N \right|\right|^2_{V^2} d\mathbb{P} dx.\]
Since the inner integral is a non-negative quantity, this is
\[\geq \int_U \mathbb{E}\left[ \left|\left|\{X_n\phi_n(x)\}_{n=1}^N\right|\right|^2_{V^2}\right] dx.\]

We consider a fixed $x \in U$. By definition of $U$, we have $|\phi_n(x)|\geq \delta$ for at least $\epsilon' N$ values of $n$. We now define new independent Gaussian random variables $Y_1, \ldots, Y_{\widetilde{N}}$ for $\widetilde{N} \geq \epsilon' N$ as follows. We start from $n=1$, and we define $Y_1$ to be the first partial sum $\sum_{n=1}^{n_1} \phi_n(x) X_n$ such that $\sum_{n=1}^{n_1} |\phi_n(x)| \geq \delta$. We then similarly define $Y_2$ to be $\sum_{n=n_1 +1}^{n_2} \phi_n(x) X_n$ for the smallest $n_2$ such that $\sum_{n = n_1+1}^{n_2} |\phi_n(x)|\geq \delta$. We continue this process, defining the $Y_i$'s to be disjoint sums of the $\phi_n(x) X_n$'s. Since $x \in U$, we will have $Y_1, \ldots, Y_{\widetilde{N}}$ for $\widetilde{N} \geq \epsilon' N$. Since the sum of independent Gaussians is distributed as a Gaussian (with variance equal to the sum of the variances), each $Y_i$ is distributed as an independent, mean zero Gaussian with variance $\geq \delta^2$.
Thus, applying Lemma \ref{lem:Qian}, we have for each $x \in U$:
\[\mathbb{E}\left[ \left|\left| \{X_n\phi_n(x)\}_{n=1}^N\right|\right|_{V^2}^2 \right] \geq \mathbb{E}\left[ \left|\left| \{Y_i\}_{i=1}^{\widetilde{N}}\right|\right|^2_{V^2} \right] \geq \delta^2 \widetilde{N} \ln \ln (\widetilde{N}) \gg \delta^2 N \ln \ln (N).\]

Therefore, we have
\begin{equation}\label{expectationbound}
\mathbb{E}\left[ \left|\left| \{X_n \phi_n(x)\}_{n=1}^N \right|\right|^2_{L^2(V^2)}\right] \gg \int_U \delta^2 N \ln \ln (N) dx \gg N \ln \ln N.
\end{equation}
We note that the constants being subsumed by the $\gg$ notation above depend on $M$.

Now, we consider the contribution to this expectation from points $\omega$ in the probability space $\Omega$ such that $\sum_{n=1}^N X_n(\omega)^2$ is much larger than $N$. We will show this contribution is small. To do this, we will upper bound the quantity $\mathbb{P}\left[\sum_{n=1}^N X_n^2 \geq kN\right]$ for each positive integer $k\geq 2$.
We rely on the following version of the Berry-Esseen theorem.

\begin{lemma}\label{lem:strongbe}(\cite{Petrov}, p. 132) Let $Z_1, \ldots, Z_N$ be independent, mean zero random variables with $\mathbb{E}[|Z_n|^{2+\gamma}] < \infty$ for all $n$ for some $0 < \gamma \leq 1$. Let $\sigma^2_n := \mathbb{E}[Z_n^2]$ and $B_N := \sum_{n=1}^N \sigma_n^2$. Then, for all $x \in \R$:
\[\left| \mathbb{P}\left[ B_N^{-\frac{1}{2}} \sum_{n=1}^N Z_n< x\right] - \Phi(x) \right| \leq \frac{A}{B_N^{1+\gamma/2} (1+|x|)^{2+\gamma}} \sum_{n=1}^N \mathbb{E}[|Z_n|^{2+\gamma}],\]
where $A$ is a constant and $\Phi(x)$ denotes the standard normal distribution function.
\end{lemma}

Now, letting $X_1, \ldots, X_N$ denote the independent, mean zero, variance one Gaussians as above, we define $Z_1, \ldots, Z_N$ by $Z_n := X_n^2-1$. Then the $Z_n$'s are independent, mean zero random variables. We note that $\mathbb{E}[Z_n^2] = \mathbb{E}[X_n^4] -1 = 2$ for each $n$. Also,
\[\mathbb{E}[|Z_n|^3] = \mathbb{E}[|X_n^6-3X_n^4+3X_n^2-1|] \leq \mathbb{E}[X_n^6] + 3\mathbb{E}[X_n^4] + 3\mathbb{E}[X_n^2]+1 = 28.\]
We will apply Lemma \ref{lem:strongbe} for $Z_1, \ldots, Z_N$, with $\gamma := 1$ and $B_N = 2N$ (since $\sigma_n^2 = 2$ for each $n$).
We observe:
\[\mathbb{P}\left[\sum_{n=1}^N X_n^2 \geq kN\right] = \mathbb{P}\left[\sum_{n=1}^N Z_n \geq (k-1)N\right] = \mathbb{P}\left[ B_N^{-\frac{1}{2}}\sum_{n=1}^N Z_n \geq 2^{-\frac{1}{2}}(k-1)N^{\frac{1}{2}}\right]\]
\[= 1 - \mathbb{P}\left[ B_N^{-\frac{1}{2}}\sum_{n=1}^N Z_n < x \right] \leq 1 - \Phi(x) + \frac{A}{B_N^{3/2}(1+|x|)^{3}} \sum_{n=1}^N \mathbb{E}\left[|Z_n|^3\right],\]
where $x:= 2^{-1/2}(k-1)N^{1/2}$.

Since $\mathbb{E}\left[|Z_n|^3\right]$ is a constant, this is
\[ \ll \int_{x}^\infty e^{-\frac{y^2}{2}} dy + \frac{1}{N^{1/2}(1+|x|)^3}.\]
Using that $x = 2^{-1/2} (k-1)N^{1/2}$, we have
\begin{equation}\label{errorterm}
\frac{1}{N^{1/2}(1+|x|)^3}\ll \frac{1}{N^2(k-1)^3}.
\end{equation}
Since $x \geq 1$ (recall that $k \geq 2$), we have
\begin{equation}\label{mainterm}
\int_{x}^\infty e^{-\frac{y^2}{2}}dy \leq \int_x^\infty y e^{-\frac{y^2}{2}} dy = e^{-\frac{x^2}{2}}= e^{-\frac{1}{4} N(k-1)^2}.
\end{equation}
Combining (\ref{errorterm}) and (\ref{mainterm}), we see that
\[\mathbb{P}\left[\sum_{n=1}^N X_n^2 \geq kN\right] \ll \frac{1}{N^2(k-1)^3} + e^{-\frac{1}{4}N(k-1)^3},\]
for each positive integer $k \geq 2$.

Now, by Lemma \ref{varRM1}, for each $\omega \in \Omega$ such that $ kN \leq \sum_{n=1}^N X_n^2(\omega) < (k+1)N$, we have that the quantity $\left|\left| \{X_n\phi_n(x)\}_{n=1}^N\right|\right|^2_{L^2(V^2)}$ evaluated at $\omega$ is $\ll (k+1)\ln^2(N)N$. Thus, the contribution to the expectation bounded in (\ref{expectationbound}) coming from such points $\omega$ for all $k\geq 2$ is upper bounded as:
\[ \ll \sum_{k=2}^{\infty} (k+1)\ln^2(N) N \left( e^{-\frac{1}{4}N(k-1)^2}+ \frac{1}{N^2(k-1)^3}\right)\]
\[ = \ln^2(N) N e^{-\frac{1}{4}N} \sum_{k=2}^\infty (k+1)\left(e^{-\frac{1}{4}N}\right)^{k^2-2k} + \frac{\ln^2(N)}{N}\sum_{k=2}^{\infty} \frac{k+1}{(k-1)^3}.\]
Both of these sums are convergent, and it is easy to see that this quantity is $o(N \ln \ln N)$.

Therefore, by (\ref{expectationbound}) and the above bounds, we have proven that there exists some point $\omega \in \Omega$ such that when we define $a_n:= X_n(\omega)$ and define $f(x) =\sum_{n=1}^{N} a_n \phi_n (x)$, we have
\[ ||S[f]||_{L^2(V^2)} \gg_{M} \sqrt{\ln\ln(N)}||f||_{L^2}. \]
Here, we have used that we can choose $\omega$ so that $||S[f]||_{L^2(V^2)}^2 \gg_M N \ln \ln (N)$ and $||f||^2_{L^2} = \sum_{n=1}^N a_n^2 \leq 2N$ simultaneously.

\section{Systems of Bounded Independent Random Variables}
\label{sec:probability}

In this section, we prove the following theorem:

\begin{thm9} Let $\{X_i\}_{i=1}^N$ be a sequence of mean zero independent random variables such that $|X_i|\leq C$ and $\mathbb{E}\left[ |X_i|^2 \right] =1$ for all $i \in [N]$. Then
\[\mathbb{E}\left[\left| \left|  \{a_i X_i\}_{i=1}^N\right|\right|_{V^2}\right] \ll_C \sqrt{\ln\ln(N)} \left(\sum_{i=1}^{N} a_i^2 \right)^{1/2} . \]
\end{thm9}

We will require the following lemmas. The first is a form of Hoeffding's inequality \cite{Hoeffding}.

\begin{lemma}\label{lem:Hoeffding}Let $\{X_i\}$ be independent random variables such that $\mathbb{P}[X_i \in [a_i,b_i]]=1$. Then
\[\mathbb{P}\left[ \left| S_n - \mathbb{E}\left[S_n\right] \right| \geq t \right] \leq 2\exp\left(  - \frac{2t^2}{\sum_{i=1}^{n}(b_i-a_i)^2} \right) \]
where $S_n = \sum_{i=1}^{n} X_i$.
\end{lemma}

\begin{lemma}\label{lem:etemadi} (Etemadi's Inequality). (See Theorem 1 in \cite{E}.) Let $X_1, X_2, \ldots, X_n$ denote independent random variables and let $a > 0$. Let $S_\ell := X_1 + \cdots + X_\ell$ denote the partial sum. Then
\[\mathbb{P} [\max_{1 \leq \ell \leq n} |S_\ell| \geq 3 a] \leq 3 \max_{1 \leq \ell \leq n} \mathbb{P}[|S_\ell| \geq a].\]
\end{lemma}

\begin{lemma}\label{lem:Rosenthal} (Rosenthal's Inequality). (See Theorem 3 in \cite{R}.) Let $2 < p < \infty$. Then there exists a constant $K_p$ depending only on $p$, so that if $X_1, \ldots, X_n$ are independent random variables with $\mathbb{E}[X_i] = 0$ for all $i$ and $\mathbb{E}[|X_i|^p] < \infty$ for all $i$, then:
\[\left(\mathbb{E}[|S_n|^p]\right)^{1/p} \leq K_p\; \max \left\{ \left( \sum_{i=1}^{n} \mathbb{E}[|X_i|^p]\right)^{1/p}, \left(\sum_{i=1}^{n} \mathbb{E}[|X_i|^2]\right)^{1/2}\right\}.\]
\end{lemma}

We also use the following consequence of Doob's inequality. For an interval $I \subseteq [n]$, we define $S_I:= \sum_{i\in I} X_i$. We also define
\[\tilde{S}_n := \max_{I \subseteq [n]} |S_I|.\] We then have:
\begin{lemma}\label{lem:Doob} For $p >1$ and independent random variables $X_1, \ldots, X_n$ with $\mathbb{E}[X_i] = 0$ for all $i$,
\[\mathbb{E}\left[ |\tilde{S}_n|^p\right] \leq 2^p \mathbb{E}\left[ \max_{1\leq \ell \leq n} \left| \sum_{i=1}^\ell X_i\right|^p\right] \leq 2^p \left(\frac{p}{p-1}\right)^p \mathbb{E}\left[ |S_n|^p\right].\]
\end{lemma}

\begin{proof}
The first inequality is a consequence of the following observation. For a subinterval $I \subseteq [n]$, we let $I_0$ be the subinterval that starts at 1 and ends just before $I$, and we let $I_1$ be the interval $I_0 \cup I$. Then $I_0$ and $I_1$ are both intervals starting at 1, and $S_{I_0} + S_{I} = S_{I_1}$. Therefore, $\max \{ |S_{I_0}|, |S_{I_1}|\} \geq \frac{1}{2} |S_I|$. The second inequality follows from Theorem 3.4 on p. 317 in \cite{Doob}.
\end{proof}

We begin by decomposing $[N]$ into a family of subintervals according to a concept of mass defined with respect to the $a_i$ values.
We  define the \emph{mass} of a subinterval $I \subseteq [N]$ as $M(I) := \sum_{n \in I} a_{n}^2$. By normalization, we may assume that $M([N])=1$. We define $I_{0,1} := [N]$ and we iteratively define $I_{k,s}$, for $1\leq s\leq 2^k$, as follows. Assuming we have already defined $I_{k-1,s}$ for all $1 \leq s \leq 2^{k-1}$, we will define $I_{k,2s-1}$ and $I_{k,2s}$, which are subintervals of $I_{k-1,s}$. $I_{k,2s-1}$ begins at the left endpoint of $I_{k-1,s}$ and extends to the right as far as possible while covering strictly less than half the mass of $I_{k-1,s}$, while $I_{k,2s}$ ends at the right endpoint of $I_{k-1,s}$ and extends to the left as far as possible while covering at most half the mass of $I_{k-1,s}$. More formally,
we define $I_{k,2s-1}$ as the maximal subinterval of $I_{k-1,s}$ which contains the left endpoint of $I_{k-1,s}$ and satisfies $M(I_{k,2s-1}) < \frac{1}{2} M(I_{k,s})$. We also define $I_{k,2s}$ as the maximal subinterval of $I_{k-1,s}$ which contains the right endpoint of $I_{k-1,s}$ and satisfies $M(I_{k,2s}) \leq \frac{1}{2} M(I_{k,s})$. We note that these subintervals are disjoint. We may express $I_{k-1,s} = I_{k,2s-1} \bigcup I_{k,2s} \bigcup i_{k,s}$,  where $i_{k,s} \in I_{k-1,s}$. In other words, $i_{k,s}$ denotes the single element which lies between $I_{k,2s-1}$ and $I_{k,2s}$ (note that such a point always exists because we have required that $I_{k,2s-1}$ contains strictly less than half of the mass of the interval). Here it is acceptable, and in many instances necessary, for some choices of the intervals in this decomposition to be empty.
By construction we have that

\begin{equation}\label{eq:mass}
M(I_{k,s}) \leq 2^{-k}.
\end{equation}

We call an interval $J \subseteq [N]$ admissible if it is an element of the decomposition given above. We denote the collection of admissible intervals by $\mathcal{A}$. We additionally refer to the subset $\{I_{k,s}| 1\leq s\leq 2^k\}$ of $\mathcal{A}$ as the admissible intervals on level $k$ and the subset $\{i_{k,s} | 1 \leq s \leq 2^{k}\}$ as the admissible points on level $k$. We note that every point in $[N]$ is an admissible point on some level. (Eventually, we have subdivided all intervals down to being single elements.)

We consider an arbitrary interval $J \subseteq [N]$. We would like to approximate $J$ by an admissible interval $\tilde{J}$ such that $J \subseteq \tilde{J}$ and $M(\tilde{J}) \leq c M(J)$, for some constant $c$. This may be impossible, however, since $J$ could span the boundary between adjacent admissible intervals for all comparable masses. To address this, we will instead approximate $J$ by the union of two admissible intervals and one point.

\begin{lemma}\label{lem:decomposition} For every $J \subseteq [N]$, ($J \neq \emptyset$) there exist $\tilde{J}_\ell, \tilde{J}_r \in \mathcal{A}$ and $i_J \in [N]$ such that $\tilde{J}:= \tilde{J}_{\ell} \cup i_J \cup \tilde{J}_r$ is an interval (i.e. $J_\ell, i_J, J_\ell$ are adjacent), $J \subseteq \tilde{J}$, and $M(\tilde{J}) \leq 2M(J)$.
\end{lemma}

\begin{proof} We consider the minimal value $k$ such that $J$ contains an admissible point on level $k$. We note that this point is unique, and we define $i_J$ to be equal to it. To see why a unique such point exists, first note that if $J$ contained at least two admissible points on level $k$, then it would also contain an admissible point between them on level $k-1$. Now we consider the subinterval $J_\ell$ consisting of elements of $J$ that lie to the left of $i_j$. Since the rightmost endpoint of this subinterval is at rightmost endpoint of an admissible interval on level $k$, it is also a rightmost endpoint of some admissible interval on every level $> k$. We define $\tilde{J}_\ell$ to be the admissible interval with this right endpoint on the highest level $k_\ell$ such that $J_\ell \subseteq \tilde{J}_\ell$. We note that the admissible interval with this right endpoint on level $k$ contains $J$, so such an interval $\tilde{J}_\ell$ must exist, and $k_\ell \geq k$.

We claim that $M(\tilde{J}_\ell) \leq 2M(J_\ell)$. To prove this, we consider the admissible interval $\tilde{J}'$ on level $k_\ell+1$ with this same right endpoint. By maximality of $k_\ell$, we must have that $J \nsubseteq \tilde{J}'$. This implies that $J$ must contain the admissible point on level $k_\ell+1$ that occurs when $\tilde{J}_\ell$ is decomposed. Therefore, $M(J_\ell) \geq \frac{1}{2} M(\tilde{J}_\ell)$.

We define the subinterval $J_r$ consisting of elements of $J$ that lie to the right of $i_j$, and we can similarly find an admissible $\tilde{J}_r$ such that $J_r \subseteq \tilde{J}_r$ and $M(\tilde{J}_r) \leq 2M(J_r)$. We then have $J \subseteq \tilde{J}:= \tilde{J}_\ell \cup i_J \cup \tilde{J}_r$ and $M(\tilde{J})\leq 2M(J)$ follows from:
\[M(\tilde{J}) = M(\tilde{J}_\ell) + M(i_J) + M(\tilde{J}_r) \leq 2(M(J_\ell) + M(i_J) + M(J_r)) = 2M(J).\]
\end{proof}

Defining $\tilde{J}_\ell$, $\tilde{J}_r$, and $i_J$ with respect to $J$ as in the lemma, we observe that:
\begin{equation}\label{SquareDecomp}
 |S_{J}|^2 \ll |\tilde{S}_{\tilde{J}_{\ell}}|^2 + |\tilde{S}_{\tilde{J}_{r}}|^2 + |S_{i_{J}}|^2.
\end{equation}
Here, $|\tilde{S}_{\tilde{J}}|$ is the maximal partial sum over all subintervals contained in $\tilde{J}$.
Also, if $\mathcal{P}$ is a partition of $[N]$, then the admissible intervals and points ($\tilde{J}_{\ell}$, $\tilde{J}_{r}$, and $i_{J}$) associated to an element $J$ of the partition will only reoccur for a bounded number of elements of the partition (i.e. a particular admissible interval/point will only appear among $\tilde{J}_\ell, \tilde{J}_r, i_J$ for a constant number of $J \in \mathcal{P}$). This is because the $J$'s in $\mathcal{P}$ are disjoint, so $i_J \in J$ for only one $J \in \mathcal{P}$, and $M(J\cap \tilde{J}_\ell) \geq \frac{1}{2} \tilde{J}_\ell$ implies $\tilde{J}_\ell$ can appear for at most two $J$'s in $\mathcal{P}$.

Now we will prove Theorem \ref{varRad}. We let $\Omega$ denote the probability space for $X_1, \ldots, X_N$ (each $\omega$ in $\Omega$ is associated to a sequence of $N$ real numbers). For each $\omega \in \Omega$, we let $\mathcal{P}_{\omega}$ denote a maximizing partition. We define $\mathcal{P}_{\omega,\ell}$ (resp. $\mathcal{P}_{\omega,r}$) to be the set of $\tilde{J}_{\ell}$ (resp. $\tilde{J}_{r}$) associated to $J \in \mathcal{P}_{\omega}$. We note that the same interval could appear as $\tilde{J}_\ell$ or $\tilde{J}_r$ for up to two different $J$'s in $\mathcal{P}_\omega$.

We fix a large constant $B$ which will be specified later. Now we split each set $\mathcal{P}_{\omega,\text{side}}$ (here $\text{side} \in \{\ell,r\}$) into two disjoint subsets $\mathcal{P}_{\omega,\text{side}}^{\text{good}}$ and $\mathcal{P}_{\omega,\text{side}}^{\text{bad}}$. We define $\mathcal{P}_{\omega,\text{side}}^{\text{good}}$ to be the set of $\tilde{J} \in \mathcal{P}_{\omega,\text{side}}$ such that
\begin{equation}
 \left|\tilde{S}_{\tilde{J}} \right|^2 \leq B M(\tilde{J}) \ln \ln (N).
\end{equation}
We then define $\mathcal{P}_{\omega,\text{side}}^{\text{bad}}$ to be the complement of $\mathcal{P}_{\omega,\text{side}}^{\text{good}}$ inside $\mathcal{P}_{\omega,\text{side}}$.

Our objective is to prove the estimate
\[ \mathbb{E} \left[  \sum_{J \in \mathcal{P}_{\omega}} |S_{J}|^2 \right] \ll  \ln \ln (N) .\]
Using (\ref{SquareDecomp}), we upper bound the left side as follows:
\[\mathbb{E} \left[  \sum_{J \in \mathcal{P}_{\omega}} |S_{J}|^2 \right]  \ll
\mathbb{E} \left[  \sum_{\tilde{J} \in \mathcal{P}_{\omega,l}^{\text{good}} } |\tilde{S}_{\tilde{J}}|^2 \right]  +
\mathbb{E} \left[  \sum_{\tilde{J} \in \mathcal{P}_{\omega,r}^{\text{good}} } |\tilde{S}_{\tilde{J}}|^2 \right]\]
\[+ \mathbb{E} \left[  \sum_{\tilde{J} \in \mathcal{P}_{\omega,l}^{\text{bad}} } |\tilde{S}_{\tilde{J}}|^2 \right] +
\mathbb{E} \left[  \sum_{\tilde{J} \in \mathcal{P}_{\omega,r}^{\text{bad}} } |\tilde{S}_{\tilde{J}}|^2 \right]
 + \mathbb{E} \left[  \sum_{J \in \mathcal{P}_{\omega}} |S_{i_J}|^2 \right].
\]

We observe that $ \sum_{\tilde{J} \in \mathcal{P}_{\omega,\text{side}}^{\text{good}}} |\tilde{S}_{\tilde{J}}|^2 \ll \left( \sum_{\tilde{J} \in \mathcal{P}_{\omega,\text{side}}} M(\tilde{J}) \right) \ln\ln (N) \ll \ln\ln (N) $. This holds because $\sum_{J \in \mathcal{P}} M(J) = 1$, and the total mass of the intervals $\tilde{J}_\ell, \tilde{J}_r, i_J$ used to cover each $J$ is at most $2M(J)$, thus $\sum_{\tilde{J} \in \mathcal{P}_{\omega, \text{side}}} M(\tilde{J}) \leq 2$.
This shows that the terms involving the good admissible intervals are easily controlled. The last term is also easily controlled as follows
\[\mathbb{E} \left[ \sum_{J \in \mathcal{P}_{\omega}} |S_{i_J}|^2 \right]  \ll \mathbb{E} \left[ \sum_{n \in [N]}  |a_n X_{n}|^2 \right] \ll 1. \]

It remains to control the terms involving the bad admissible intervals. The argument is essentially the same for both the sums over $\mathcal{P}_{\omega,l}^{\text{bad}}$ and $\mathcal{P}_{\omega,r}^{\text{bad}}$, so we will work with the quantity $\mathbb{E} \left[  \sum_{\tilde{J} \in \mathcal{P}_{\omega,\text{side}}^{\text{bad}} } |\tilde{S}_{\tilde{J}}|^2 \right]$ in what follows.

We now partition $\mathcal{P}_{\omega,\text{side}}^{\text{bad}}$ into two disjoint sets $\mathcal{P}_{\omega,\text{side}}^{\text{bad},1}$ and $\mathcal{P}_{\omega,\text{side}}^{\text{bad},2}$. The set $\mathcal{P}_{\omega,\text{side}}^{\text{bad},1}$ consists of intervals $I_{k,s} \in \mathcal{P}_{\omega,\text{side}}^{\text{bad}}$ such that $|I_{k,s}| \leq 2^{-k/2} N$ and $\mathcal{P}_{\omega,\text{side}}^{\text{bad},2}$ contains the complement set.
For each $k$, we define $T_{k} \subseteq \{I_{k,s} : 1\leq s \leq 2^{k}\} $ as the collection of all intervals $I_{k,s}$ satisfying $|I_{k,s}| \geq 2^{-k/2} N$. Clearly, $|T_{k}| \leq 2^{k/2}$ for each $k$.
We then have:
\[\mathbb{E} \left[  \sum_{\tilde{J} \in \mathcal{P}_{\omega,\text{side}}^{\text{bad},2} } |\tilde{S}_{\tilde{J}}|^2 \right] \ll
\mathbb{E} \left[  \sum_{k=1}^{\infty} \sum_{\tilde{J} \in T_{k}} |\tilde{S}_{\tilde{J}}|^2 \right] =  \sum_{k=1}^{\infty} \sum_{\tilde{J} \in T_{k}} \mathbb{E} \left[ |\tilde{S}_{\tilde{J}}|^2 \right].
\]

Using (\ref{eq:mass}) and the fact that $ \mathbb{E} \left[ |\tilde{S}_{\tilde{J}}|^2\right] \ll \mathbb{E} \left[ |S_{\tilde{J}}|^2 \right]$ (by Lemma \ref{lem:Doob}), we have
\[ \sum_{k=1}^{\infty} \sum_{\tilde{J} \in T_{k}} \mathbb{E} \left[ |\tilde{S}_{\tilde{J}}|^2 \right] \ll \sum_{k=1}^{\infty}  \sum_{\tilde{J} \in T_{k}} \mathbb{E} \left[ |S_{\tilde{J}}|^2 \right] \ll \sum_{k=1}^{\infty} 2^{k/2} 2^{-k} \ll 1.
\]

It now suffices to bound the more difficult term $\mathbb{E} \left[  \sum_{\tilde{J} \in \mathcal{P}_{\omega,\text{side}}^{\text{bad},1} } |\tilde{S}_{\tilde{J}}|^2 \right]$.

Now $|I_{k,s}| \leq 2^{-k/2}N$ if $I_{k,s} \in \mathcal{P}_{\omega,\text{side}}^{\text{bad},1}$. For a fixed interval $J$, we let $B(J) \subseteq \Omega$ denote the event that the $|\tilde{S}_{J}(\omega)|^2$ is bad. In other words, $\omega \in B(J)$ if $\left| \tilde{S}_{J}(\omega) \right|^2 \geq B M(J) \ln \ln (N)$. We let $T_{k}^{c}$ denote the complement of $T_{k}$.  We now have that
\[\mathbb{E} \left[  \sum_{\tilde{J} \in \mathcal{P}_{\omega,\text{side}}^{\text{bad},1} } |\tilde{S}_{\tilde{J}}|^2 \right] \ll
 \sum_{k=1}^{2 \ln(N) } \sum_{\tilde{J} \in T_{k}^{c}} \mathbb{E} \left[ |1_{B(\tilde{J})} \tilde{S}_{\tilde{J}}|^2 \right] .\]
Here we have restricted the the summation of $k$ to the range $1\leq k \leq 2 \ln(N)$ using the fact that $1 \leq |I_{k,s}| \leq 2^{-k/2}N$ implies $k \leq 2 \ln (N)$.

We let $\gamma > 0$ denote a positive value to be specified later.
Letting $2p:=2+\gamma$ and applying Lemma \ref{lem:Rosenthal} (Rosenthal's inequality) we have that
\[ \left(\mathbb{E} \left[ |  S_{\tilde{J}} |^{2p} \right] \right)^{1/p} =\left(\mathbb{E}\left[ \left| S_{\tilde{J}} \right|^{2+\gamma} \right]\right)^{\frac{2}{2+\gamma} } \ll \left( \mathbb{E}\left[ \left| \sum_{n\in \tilde{J}} a_n X_{n} \right|^{2+\gamma} \right]\right)^{\frac{2}{2+\gamma} } \]
\begin{equation}\label{holder1} \ll \max \left\{ \left(\sum_{n\in \tilde{J}} |a_n|^{2+\gamma} \mathbb{E}\left[ |X_{i}| ^{2+\gamma}\right]\right)^{\frac{2}{2+\gamma}} , \left( \sum_{n\in \tilde{J}} |a_n|^2\right)\right\}\ll  M(\tilde{J}).
\end{equation}

The last inequality follows from the fact that the $\ell^{2}$ norm is greater than the $\ell^{2+\gamma}$ norm and $\mathbb{E}\left[|X_i|^{2+\gamma}\right] \leq C^{2+\gamma}$.

%

We let $s:= |\tilde{J}|$, and we let $S_{\tilde{J}, \ell}$ denote the sum of $a_iX_i$ for the first $\ell$ indices $i$ in $\tilde{J}$. By definition of the event $B(\tilde{J})$, we have:
\[\mathbb{E} \left[ 1_{B(\tilde{J})}\right] = \mathbb{P}\left[\left| \tilde{S}_{\tilde{J}} \right|^2 \geq B M(\tilde{J}) \ln \ln (N) \right] \leq \mathbb{P}\left[\max_{1\leq \ell \leq s} \left| S_{\tilde{J},\ell}\right|^2 \geq \frac{B}{2} M(\tilde{J}) \ln \ln (N)\right].\]
By Lemma \ref{lem:etemadi}, this is
\[ \ll \max_{1 \leq \ell \leq s} \mathbb{P}\left[ \left| S_{\tilde{J},\ell} \right|^2 \geq \frac{B}{6} M(\tilde{J})\ln \ln (N)\right].\]

By Lemma \ref{lem:Hoeffding}, this is:
\[\ll \exp\left(  - \frac{ B M(\tilde{J}) \ln \ln (N) }{ 3C^2 M(\tilde{J}) } \right) = \exp \left( - \frac{B \ln \ln (N)}{3C^2}\right).\]
By setting the value of $B$ to be sufficiently large with respect to the constant $C$ (i.e. $B > 12 C^2$), we have:
\begin{equation}
    \mathbb{E} \left[ 1_{B(\tilde{J})}\right]\ll \ln^{-4}(N).
\end{equation}

We now define $q$ as a function of $p$ so that $\frac{1}{p} + \frac{1}{q} = 1$, i.e. $q = \frac{p}{p-1}$. We then set $\gamma$ such that
\begin{equation}\label{holder2}
\left( \mathbb{E} \left[ 1_{B(\tilde{J})} \right] \right)^{1/q} \ll \ln^{-2}(N)
\end{equation}
for all $\tilde{J}$. (Recall that $p:= \frac{2+\gamma}{2}$.)
We now apply H\"{o}lder's inequality with $p$ and $q$ to obtain:
\[\sum_{k=1}^{2\ln (N) } \sum_{\tilde{J} \in T_k^c} \mathbb{E} \left[ \left|1_{B(\tilde{J})}\tilde{S}_{\tilde{J}}^2\right|\right]
\leq \sum_{k=1}^{2 \ln (N)} \sum_{\tilde{J} \in T_k^c} \left( \mathbb{E} \left[ \left| 1_{B(\tilde{J})}\right|^q \right] \right)^{\frac{1}{q}} \left( \mathbb{E} \left[ \left| \tilde{S}_{\tilde{J}}\right|^{2p}\right]\right)^{\frac{1}{p}}.\]
Using (\ref{holder1}), (\ref{holder2}) and Lemma \ref{lem:Doob}, we see this is:
\[\ll \sum_{k=1}^{2\ln (N)} \sum_{\tilde{J} \in T_k^c} \ln^{-2}(N) M(\tilde{J}) \ll \sum_{k=1}^{2 \ln (N)} \ln^{-2}(N) \ll \frac{1}{\ln (N)}.\]

This completes the proof.

\section{Random Permutations}

In this section, we will use probabilistic techniques to prove the following theorem:

\begin{thm7} Let $\{ \phi_n \}_{n=1}^N$ be an orthonormal system such that $|\phi_n(x)| =1$ for all $n$ and all $x \in \T$, and $\{ a_n \}_{n=1}^N$ a choice of (complex) coefficients. Then there exists a permutation $\pi:[N] \rightarrow [N]$ such that
\[
 \left|\left|\{ a_{\pi(n)}\phi_{\pi(n)}\}_{n=1}^N\right|\right|_{L^2(V^2)} \ll \sqrt{\ln\ln(N)} \left(\sum_{n=1}^{N} |a_n|^2 \right)^{1/2}  \]
\end{thm7}

\begin{proof}
We assume without loss of generality that $\sum_{n=1}^N |a_n|^2 = 1$. Then, for each $a_n$, there exists some non-negative integer $j$ such that $2^{-j-1} < |a_n|^2 \leq 2^{-j}$. For each fixed $j$, we let $A_j$ denote the set of $n \in [N]$ such that $2^{-j-1} < |a_n|^2 \leq 2^{-j}$.
We define $A^* \subseteq [N]$ as $A^*:= \bigcup_{j=\lceil 2\ln N \rceil}^\infty A_j$. We also define
\[b_n = \left\{
    \begin{array}{ll}
      a_n, & \hbox{$n \in A^*$} \\
      0, & \hbox{$n \notin A^*$.}
    \end{array}
  \right.\]
We then observe, for any permutation $\pi:[N]\rightarrow [N]$ and any $x \in \T$,
\[\left|\left| \{b_{\pi(n)} \phi_{\pi(n)}(x)\}_{n=1}^N\right|\right|_{V^2} \ll \sum_{n=1}^N \left| b_n \phi_n(x)\right| \ll \frac{1}{N}\cdot N \ll 1.\]
Applying the triangle inequality for the $||\cdot ||_{V^2}$ norm, this allows us to ignore the contribution of all terms $a_n$ where $n \in A^*$.

We consider the class of permutations $\pi: [N] \rightarrow [N]$ such that $\pi^{-1}(A_j)$ is an interval for each $j$. In other words, these are permutations which group the elements of each $A_j$ together. We allow arbitrary orderings within each group and an arbitrary ordering of the groups. For a fixed permutation $\pi$, we let $B_j$ denote the preimage of $A_j$ under $\pi$ (so $B_j$ is an interval). We will refer to the intervals $B_j$ as ``blocks".
From this point onward, we will only consider permutations belonging to this class, and we will only consider the contribution of terms for $A_1$ up to $A_{\lfloor 2 \ln (N) \rfloor}$. We let $N' := |A_1| + \cdots + |A_{\lfloor 2 \ln(N) \rfloor}|$. For notational convenience, we assume that $\pi$ maps $[N']$ bijectively to $\bigcup_{i=1}^{\lfloor 2\ln (N)\rfloor} A_j$. (This is without loss of generality, since we have seen that we can treat the set $A^*$ separately.)

For each fixed permutation $\pi:[N]\rightarrow [N]$ in this class and each fixed $x \in \T$, we consider the quantity
\begin{equation}\label{allpartitions}
\left| \left| \{a_{\pi(n)} \phi_{\pi(n)}(x)\}_{n=1}^{N'}\right|\right|^2_{V^2} = \sum_{I \in \mathcal{P}} \left|\sum_{n \in I} a_{\pi(n)} \phi_{\pi(n)}(x)\right|^2,
\end{equation}
where $\mathcal{P}$ denotes the maximizing partition of $[N']$.

We now define two additional operators, $V^2_L$ and $V^2_S$. The value of $\left|\left| \{a_{\pi(n)} \phi_{\pi(n)}(x)\}_{n=1}^{N'}\right|\right|^2_{V^2_L}$ is defined as
\[\left|\left| \{a_{\pi(n)} \phi_{\pi(n)}(x)\}_{n=1}^{N'} \right| \right|^2_{V^2_L} := \sum_{I \in \mathcal{P}_L} \left|\sum_{n\in I} a_{\pi(n)} \phi_{\pi(n)}(x) \right|^2,\]
where $\mathcal{P}_L$ is the maximizing partition among the subset of partitions of $[N']$ that use only intervals which are unions of the $B_j$'s.

The value of $\left|\left| \{a_{\pi(n)} \phi_{\pi(n)}(x)\}_{n=1}^{N'}\right|\right|^2_{V^2_S}$ is defined as
\[\left|\left| \{a_{\pi(n)} \phi_{\pi(n)}(x)\}_{n=1}^{N'}\right|\right|^2_{V^2_S} := \sum_{I \in \mathcal{P}_S} \left|\sum_{n\in I} a_{\pi(n)} \phi_{\pi(n)}(x) \right|^2,\]
where $\mathcal{P}_S$ is the maximizing partition among the subset of partitions of $[N']$ that use only intervals $I$ that are contained in some $B_j$. This can be alternatively described as taking that maximizing partition of each $B_j$ and then taking a union of these to form $\mathcal{P}_S$.

We now claim:
\begin{equation}\label{longshort}
    \left| \left| \{a_{\pi(n)} \phi_{\pi(n)}(x)\}_{n=1}^{N'}\right|\right|^2_{V^2} \ll \left|\left| \{a_{\pi(n)} \phi_{\pi(n)}(x)\}_{n=1}^{N'}\right|\right|^2_{V^2_L} + \left|\left| \{a_{\pi(n)} \phi_{\pi(n)}(x)\}_{n=1}^{N'}\right|\right|^2_{V^2_S}.
\end{equation}

To see this, consider the maximizing partition $\mathcal{P}$ in (\ref{allpartitions}). Each $I \in \mathcal{P}$ can be expressed as the union of three disjoint intervals, $I_{S_\ell}$, $I_L$, and $I_{S_r}$, where $I_{S_\ell}$ and $I_{S_r}$ are each contained in some $B_i$, and $I_L$ is a union of $B_i$'s. More precisely, $I_L$ is the union of all the intervals $B_j$ that are contained in $I$, $I_{S_\ell}$ goes from the left endpoint of $I$ until the left endpoint of $I_L$, and $I_{S_r}$ goes from the right endpoint of $I_L$ until the right endpoint of $I$. By construction, each of $I_{S_\ell}$ and $I_{S_r}$ is contained in some $B_j$. (Some of $I_L, I_{S_r}, I_{S_\ell}$ may be empty.) Thus,
\[\left|\sum_{n \in I} a_{\pi(n)} \phi_{\pi(n)}(x)\right|^2 \ll \left| \sum_{n\in I_L} a_{\pi(n)} \phi_{\pi(n)}(x)\right|^2 + \left|\sum_{n \in I_{S_\ell}} a_{\pi(n)} \phi_{\pi(n)}(x)\right|^2 + \left|\sum_{n \in I_{S_r}} a_{\pi(n)} \phi_{\pi(n)}(x)\right|^2.\]

Now, if we consider the set of intervals $I_L$ corresponding to $I \in \mathcal{P}$, we get a disjoint set of intervals that can occur as part of a partition considered by the operator $V^2_L$. Similarly, if we consider the set of intervals $I_{S_\ell}, I_{S_r}$ corresponding to $I \in \mathcal{P}$, we get a disjoint set of intervals that can occur as part of a partition considered by the operator $V^2_S$. Therefore,
\[\sum_{I \in \mathcal{P}} \left|\sum_{n \in I} a_{\pi(n)} \phi_{\pi(n)}(x)\right|^2 \ll \sum_{I \in \mathcal{P}_L} \left|\sum_{n \in I} a_{\pi(n)} \phi_{\pi(n)}(x)\right|^2 + \sum_{I \in \mathcal{P}_S} \left|\sum_{n \in I} a_{\pi(n)} \phi_{\pi(n)}(x)\right|^2.\]
The inequality (\ref{longshort}) then follows.

We first bound the contribution of the $V^2_L$ operator. For each $B_j$, we define the function $f_j: \T \rightarrow \mathbb{C}$ as:
\begin{equation}\label{blockfunction}
f_j(x):= \sum_{n \in B_j} a_{\pi(n)} \phi_{\pi(n)}(x).
\end{equation}
Since the sets $B_j$ are disjoint, we note that the functions $f_j$ are orthogonal to each other, but they may not be uniformly bounded. We need to show that there exists a permutation $\sigma: [\lfloor 2 \ln (N) \rfloor]\rightarrow [\lfloor 2\ln(N) \rfloor]$ of the $f_j$ values such that
\begin{equation}\label{blockgoal}
\left|\left| \{f_{\sigma(j)}(x)\}_{j=1}^{\lfloor 2\ln(N)\rfloor}\right|\right|_{L^2(V^2)} \ll \sqrt{\ln \ln(N)} \left( \sum_{n=1}^N |a_n|^2\right)^{1/2}.
\end{equation}
This would imply that there is some ordering of the blocks for which the contribution of the $V^2_L$ operator is suitably bounded.

To show (\ref{blockgoal}), we will use the following inequality of Garsia for real numbers:
\begin{lemma}\label{lem:garsia}(See Theorem 3.6.15 in \cite{Garsia}.) Let $x_1, \ldots, x_M \in \mathbb{R}$. We consider choosing a permutation $\psi$ of $[M]$ uniformly at random. Then:
\[\mathbb{E}\left[ \max_{1 \leq k \leq M} \left( x_{\psi(1)} + \cdots + x_{\psi(k)}\right)^2\right] \ll \left( \sum_{k=1}^M x_k\right)^2 + \sum_{k=1}^M x_k^2.\]
\end{lemma}

We derive the following corollary:
\begin{corollary}\label{cor:garsia} Let $x_1, \ldots, x_M \in \mathbb{R}$. Let $L$ be a positive integer, $1 \leq L \leq M$. Let $\mathcal{P}$ denote the partition of $[M]$ into intervals of size $L$ (starting with $[L]$), except that the last interval may be of smaller size (when $L$ does not divide $M$). We consider choosing a permutation $\psi$ of $[M]$ uniformly at random. Then:
\[\mathbb{E} \left[ \sum_{I \in \mathcal{P}} \max_{I' \subseteq I} \left( \sum_{j\in I'} x_{\psi(j)}\right)^2 \right] \ll {M-1 \choose L-1}^{-1} \left(\sum_{\stackrel{S\subseteq [M]}{|S|=L}} \left(\sum_{j \in S} x_j\right)^2+\sum_{j \in S} x_j^2\right). \]
\end{corollary}
We note here that $S$ ranges over all subsets of $[M]$ of size $L$.

\begin{proof} By linearity of expectation, we first observe:
\[\mathbb{E} \left[ \sum_{I \in \mathcal{P}} \max_{I' \subseteq I} \left( \sum_{j\in I'} x_{\psi(j)}\right)^2 \right] = \sum_{I \in \mathcal{P}} \mathbb{E} \left[ \max_{I' \subseteq I} \left(\sum_{j \in I'} x_{\psi(j)}\right)^2\right].\]
This quantity is then
\[\ll \frac{M}{L} \; \mathbb{E}\left[ \max_{I' \subset I} \left(\sum_{j \in I'} x_{\psi(j)}\right)^2\right],\]
where $I$ is any fixed interval of size $L$ (without loss of generality, we may take $I$ to be $[L]$).

For any subset $S \subseteq [M]$ of size $L$, the probability that $\psi$ maps $I$ to $S$ is ${M \choose L}^{-1}$. Conditioned on this event, the action of $\psi$ on $I$ acts as random permutation of the values $x_j$ for $j \in S$. Applying Lemma \ref{lem:garsia}, we then have the expectation (still conditioned on $\psi$ mapping $I$ to $S$) is $\ll \left(\sum_{j \in S} x_j\right)^2 + \sum_{j \in S} x_j^2$. (Note that the maximum over all subintervals $I'$ of $I$ is bounded by a constant times the maximum over subintervals starting at the left endpoint of $I$, as in the lemma.) Thus,
\[\mathbb{E}\left[ \max_{I' \subset I} \left(\sum_{j \in I'} x_{\psi(j)}\right)^2\right] \ll {M \choose L}^{-1} \sum_{\stackrel{S \subseteq [M]}{|S|=L}} \left( \left( \sum_{j \in S} x_j\right)^2 + \sum_{j \in S} x_j^2\right).\]
Since $\frac{M}{L} {M \choose L}^{-1} = {M-1 \choose L-1}^{-1}$, the corollary follows.
\end{proof}

We now decompose $[\lfloor 2\ln (N) \rfloor]$ into a family of dyadic intervals. More precisely, we consider all dyadic intervals of the form
\[((c-1)2^\ell, c2^\ell], \; \ell \in \{0,1, \ldots, \lceil  \ln (2\ln N) \rceil\}, \; c \in \left\{1, \ldots, 2^{\lceil \ln \ln (N)+\ln 2 \rceil -\ell}\right\}\]
(Some of these intervals may go beyond $M := \lfloor 2\ln (N) \rfloor$. For these, we consider their intersection with $[M]$.) The exponent $\ell$ of an interval here defines its ``level". In other words, we say an interval $((c-1)2^\ell, c2^\ell]$ is on level $\ell$.
We let $\mathcal{F}$ denote the set of all intervals of this form.

We then have that for \emph{any} interval $I' \subseteq [M]$, there are (at most) two adjacent intervals $I_l, I_r \in \mathcal{F}$ such that $I' \subseteq I_l \cup I_r$, and $|I_l \cup I_r|\leq 4 |I'|$ (when only one interval is needed, one of $I_l, I_r$ can be substituted by $\emptyset$). To see this, consider the smallest positive integer $k$ such that $|I'| < 2^k$. Then either $I'$ is contained in some dyadic interval of length $2^k$, or it contains exactly one right endpoint of such an interval. We then take $I_l$ to the be interval on level $k$ with this right endpoint, and take $I_r$ to be the next interval (with this as its open left endpoint).

This implies the following upper bound for each permutation $\sigma$ and each $x \in \T$:
\begin{equation}\label{dyadicsum}
\left|\left| \{f_{\sigma(j)}(x)\}_{j=1}^{\lfloor 2\ln(N)\rfloor} \right|\right|^2_{V^2} \ll \sum_{I \in \mathcal{F}} \max_{I' \subseteq I} \left| \sum_{j \in I'} f_{\sigma(j)}(x)\right|^2.
\end{equation}
This holds because for each interval $J$ in the maximizing partition, $J \subseteq I_l \cup I_r$ for some $I_r, I_l \in \mathcal{F}$ with $|I| < 4|I_l \cup I_r|$. Each $I \in \mathcal{F}$ will correspond to at most a constant number of $J$'s (it can only be $I_l$ for one $J$ when $I_r$ is non-empty, $I_r$ for one $J$ when $I_l$ is non-empty, and it can contain at most 3 corresponding $J'$s), and this constant factor is absorbed by the $\ll$ notation.

We consider choosing $\sigma$ uniformly at random. We observe by Fubini's theorem:
\[\mathbb{E} \left[ \int_{\T} \left|\left|\{f_{\sigma(j)}(x)\}_{j=1}^{\lfloor 2\ln(N) \rfloor} \right|\right|^2_{V^2} dx\right] = \int_{\T} \mathbb{E}\left[\left|\left|\{f_{\sigma(j)}(x)\}_{j=1}^{\lfloor 2\ln(N) \rfloor} \right|\right|^2_{V^2}\right] dx. \]
Using the triangle inequality for the $\left|\left| \cdot \right|\right|_{V^2}$ norm and linearity of expectation, we can split each $f_j(x)$ into real and imaginary parts, $f_j(x) = f^r_j(x) + i f^i_j(x)$, where $f^r_j$ and $f^i_j$ are both real valued. We then have:
\[\ll \int_{\T} \mathbb{E}\left[\left|\left|\{f^r_{\sigma(j)}(x)\}_{j=1}^{\lfloor 2\ln(N) \rfloor} \right|\right|^2_{V^2}\right] dx + \int_{\T} \mathbb{E}\left[\left|\left|\{f^i_{\sigma(j)}(x)\}_{j=1}^{\lfloor 2\ln(N) \rfloor} \right|\right|^2_{V^2}\right] dx.\]

For each $\ell$ from 0 to $\lceil \ln (2 \ln N)\rceil$, we let $\mathcal{F}_\ell$ denote the intervals in $\mathcal{F}$ on level $\ell$. On each level, these intervals are disjoint. Applying (\ref{dyadicsum}) to the quantity above for $f^r$ (the argument for $f^i$ is identical), we can express the result as:
\[\int_{\T} \mathbb{E}\left[\left|\left|\{f^r_{\sigma(j)}(x)\}_{j=1}^{\lfloor 2\ln(N) \rfloor} \right|\right|^2_{V^2}\right] dx \ll \int_{\T} \mathbb{E} \left[ \sum_{\ell=0}^{\lceil \ln (2\ln N) \rceil} \sum_{I \in \mathcal{F}_\ell} \max_{I' \subseteq I} \left| \sum_{j \in I'} f^r_{\sigma(j)}(x)\right|^2\right] dx.\]
By linearity of expectation, this is:
\[ = \int_{\T} \sum_{\ell=0}^{\lceil \ln (2 \ln N) \rceil} \mathbb{E}\left[\sum_{I \in \mathcal{F}_\ell} \max_{I' \subseteq I} \left| \sum_{j \in I'} f^r_{\sigma(j)}(x)\right|^2\right] dx.\]

Now, for each $\ell$, we apply Corollary \ref{cor:garsia} to the dyadic intervals on level $\ell$. As a result, we see that the above quantity is
\begin{equation}\label{aftercor}
\ll \sum_{\ell=0}^{\lceil \ln (2 \ln N) \rceil} {\lfloor 2\ln(N)\rfloor -1 \choose 2^{\ell}-1}^{-1}  \sum_{\stackrel{S \subseteq [\lfloor 2\ln(N)\rfloor]}{|S| = 2^\ell}} \left( \int_{\T} \left(\sum_{j \in S} f^r_j(x)\right)^2 dx + \sum_{j \in S} \int_{\T} f^r_j(x)^2 dx\right).
\end{equation}

Combining this with the same result for the imaginary parts, we have:
\[\int_{\T} \mathbb{E}\left[\left|\left|\{f_{\sigma(j)}(x)\}_{j=1}^{\lfloor 2\ln(N) \rfloor} \right|\right|^2_{V^2}\right] dx  \ll
\sum_{\ell=0}^{\lceil \ln (2 \ln N) \rceil} {\lfloor 2\ln(N)\rfloor -1 \choose 2^{\ell}-1}^{-1} \times \]
\begin{equation}\label{backtogether}
\sum_{\stackrel{S \subseteq [\lfloor 2\ln(N)\rfloor]}{|S| = 2^\ell}} \left( \int_{\T} \left(\sum_{j \in S} f^r_j(x)\right)^2 + \left(\sum_{j\in S} f^i_j(x)\right)^2dx + \sum_{j \in S} \int_{\T} f^r_j(x)^2 + f^i_j(x)^2 dx\right)
\end{equation}

We consider the quantity
\[\int_{\T} \left(\sum_{j \in S} f^r_j(x)\right)^2 + \left(\sum_{j\in S} f^i_j(x)\right)^2dx = \int_{\T} \sum_{j, j' \in S} f^r_j(x) f^r_{j'}(x) + f^i_{j}(x)f^i_{j'}(x) dx.\]
When $j \neq j'$,
\[\int_{\T} f^r_j(x)f^r_{j'}(x) + f^i_j(x)f^i_{j'}(x) dx = 0,\]
since $f_j$ and $f_{j'}$ are orthogonal, and this is the real part of $\int_{\T} f_j(x)\overline{f_{j'}(x)} dx$.
Thus,
\[\int_{\T} \left(\sum_{j \in S} f^r_j(x)\right)^2 + \left(\sum_{j\in S} f^i_j(x)\right)^2dx \ll \sum_{j \in S} \int_{\T} f^r_j(x)^2 + f^i_j(x)^2 dx.\]

We then have:
\[\mathbb{E} \left[ \int_{\T} \left|\left|\{f_{\sigma(j)}(x)\}_{j=1}^{\lfloor 2\ln(N) \rfloor} \right|\right|^2_{V^2} dx\right] \ll \sum_{\ell=0}^{\lceil \ln (2 \ln N) \rceil} {\lfloor 2\ln(N)\rfloor -1 \choose 2^{\ell}-1}^{-1}  \sum_{\stackrel{S \subseteq [\lfloor 2\ln(N)\rfloor]}{|S| = 2^\ell}} \sum_{j \in S} \int_{\T} |f_j(x)|^2 dx.\]
By Parseval's identity, $\int_{\T} |f_j(x)|^2 dx = \sum_{n \in A_j} |a_n|^2$.
Since each $j$ occurs in exactly ${\lfloor 2\ln(N)\rfloor -1 \choose 2^{\ell}-1}$ sets of size $2^{\ell}$ for each $\ell$, the above quantity is:
\[ \ll \ln \ln (N) \sum_{n=1}^N |a_n|^2.\]
This implies that there exists some permutation $\sigma$ such that
\[\int_{\T} \left|\left|\{f_{\sigma(j)}(x)\}_{j=1}^{\lfloor 2\ln(N) \rfloor} \right|\right|^2_{V^2} dx \ll \ln \ln (N) \sum_{n=1}^N |a_n|^2.\]
Taking a square root of both sides of this establishes (\ref{blockgoal}), as desired. This concludes our analysis of the $V^2_L$ operator.

We now bound the contribution of the $V^2_S$ operator.
\begin{lemma}\label{lem:shortop} For some $\pi$ in our class of permutations,
\[\int_{\T} \left|\left|  \{a_{\pi(n)} \phi_{\pi(n)}(x)\}_{n=1}^{N'}\right|\right|^2_{V^2_S} dx \ll \ln \ln (N) \sum_{n=1}^N |a_n|^2.\]
\end{lemma}

\begin{proof}
We first observe that it suffices to prove the following inequality for each $A_j$. We let $\Pi_j$ denote the set of permutations of $A_j$, i.e. each $\pi_j \in \Pi_j$ is a bijective map from $[|A_j|] \rightarrow A_j$. We consider choosing such a permutation uniformly at random. Then if we have
\begin{equation}\label{expect}
\mathop{\mathbb{E}}_{\pi_j \in \Pi_j} \left[ \int_{\T} \left|\left| \{ a_{\pi_j(n)} \phi_{\pi_j(n)}(x)\}_{n=1}^{|A_j|}\right|\right|^2_{V^2} dx\right] \ll \ln \ln (N) \sum_{n \in A_j} |a_n|^2
\end{equation}
for each $j$, this means that there exists a permutation $\pi_j$ of each $A_j$ satisfying
\[\int_{\T} \left|\left| \{ a_{\pi_j(n)} \phi_{\pi_j(n)}(x)\}_{n=1}^{|A_j|}\right|\right|^2_{V^2} dx \ll \ln \ln (N) \sum_{n \in A_j} |a_n|^2,\]
and these permutations can be put together to form a permutation $\pi$ as required for Lemma \ref{lem:shortop}. We note that it does not matter how we concatenate the $\pi_j$'s: by definition of the $V^2_S$ operator, it only matters how each $A_j$ is permutated, not the order the $A_j$'s are placed in.

We now fix a $j$ and we will prove (\ref{expect}). By Fubini's theorem, we can interchange the order of the integral and the expectation and instead work with the quantity
\[\int_{\T} \mathop{\mathbb{E}}_{\pi_j \in \Pi_j} \left[ \left|\left| \{ a_{\pi_j(n)} \phi_{\pi_j(n)}(x)\}_{n=1}^{|A_j|}\right|\right|^2_{V^2}\right] dx.\]

For each fixed $x$, we define the set of complex numbers $\mathcal{C}$ to be the set of values $a_n \phi_n(x)$ for $n \in A_j$. Then, these complex numbers $c \in \mathcal{C}$ all satisfy $2^{-j-1} < |c|^2 \leq 2^{-j}$ (recall that $|\phi_n(x)| = 1$). We let $N_j:= |A_j|$, and we let random variables $Z_1, \ldots, Z_{N_j}$ denote random samples from $\mathcal{C}$ taken \emph{without} replacement. We then see that it suffices to show:
\begin{equation}\label{woreplace}
\mathbb{E}\left[ \left|\left| \{Z_n\}_{n=1}^{N_j}\right|\right|_{V^2}^2\right] \ll \ln \ln (N) \sum_{c \in \mathcal{C}} |c|^2 + \left|\sum_{c \in \mathcal{C}}c\right|^2.
\end{equation}

To show this, we will need the following lemma:
\begin{lemma}\label{lem:complicated} Let $X_1, \ldots, X_{N_j}$ denote uniformly random samples from $\mathcal{C}$ \textbf{with} replacement. For each $k$ from 1 to $N_j$, we let $S_k:= \sum_{i =1}^k X_i$. For a subinterval $I \subseteq [N_j]$, we let $S_I := \sum_{i \in I} X_i$. Then for any $k$ and any $p >2$:
\[\mathbb{E}\left[ \max_{I \subseteq [k]} |S_I - \mathbb{E}[S_I]|^p\right] \ll C^{p} k^{\frac{p}{2}} p^{\frac{p}{2}} 2^{-jp/2},\]
where $C$ is a positive constant.
\end{lemma}

\begin{proof} We rely on Hoeffding's inequality \cite{Hoeffding}, which implies that
\begin{equation}\label{hoeffding}
\mathbb{P}\left[ \max_{I \subseteq [k]} \left|Re[S_I] - \mathbb{E}[Re[S_I]]\right| > t\right] \ll exp\left( \frac{-ct^2}{k 2^{-j}}\right),
\end{equation}
for some positive constant $c$, where $Re[S_I]$ denotes the real part of $S_I$. (More precisely, Hoeffding's inequality is applied with the maximum over $S_m$ for $1 \leq m \leq k$. However, moving to a maximum over arbitrary subintervals only results in a change of the constant $c$.) The same holds analogously for the imaginary part of $S_I$.

We note that
\begin{equation}\label{integral}
\mathbb{E}\left[ \max_{I \subseteq [k]} \left|Re[S_I] - \mathbb{E}[Re[S_I]]\right|^p\right] = p \int_{0}^\infty t^{p-1}\mathbb{P}\left[ \max_{I \subseteq [k]} \left|Re[S_I] - \mathbb{E}[Re[S_I]]\right| > t\right] dt.
\end{equation}
Applying (\ref{hoeffding}), this is
\[\ll p \int_{0}^\infty t^{p-1} exp\left( \frac{-ct^2}{k 2^{-j}}\right) dt.\]
We now perform the change of variable $t = \lambda^{\frac{1}{p}}$, so $dt = \frac{1}{p} \lambda^{\frac{1}{p}-1} d\lambda$. We obtain:
\[ = \int_{0}^\infty exp\left( \frac{-c \lambda^{2/p}}{k 2^{-j}}\right) d\lambda.\]

We recall that $\Gamma(z):= \int_{0}^\infty t^{z-1} e^{-t} dt$. Performing the change of variable $t = s^{\frac{2}{p}}$, we have
\[\Gamma(z):= \frac{2}{p} \int_0^\infty s^{\frac{2}{p}-1} s^{\frac{2}{p} (z-1)} e^{-s^{2/p}} ds = \frac{2}{p} \int_0^\infty s^{\frac{2}{p} z -1}e^{-s^{2/p}} ds.\]
We now see that
\[\int_{0}^\infty e^{-t^{\frac{2}{p}}}dt = \frac{p}{2}\; \Gamma \left(\frac{p}{2}\right).\]

We then set $s:= \left(\frac{c}{k 2^{-j}}\right)^{p/2} \lambda$, and we have:
\[\int_{0}^\infty exp\left( \frac{-c \lambda^{2/p}}{k 2^{-j}}\right) d\lambda = \left(\frac{c}{k 2^{-j}}\right)^{-p/2} \int_{0}^\infty e^{-s^{\frac{2}{p}}} ds = \left(\frac{c}{k 2^{-j}}\right)^{-p/2} \frac{p}{2} \; \Gamma \left(\frac{p}{2}\right) .\]

This yields
\[\mathbb{E}\left[ \max_{I \subseteq [k]} \left|Re[S_I] - \mathbb{E}[Re[S_I]]\right|^p\right] \ll \frac{p}{2}  k^{p/2} c^{-p/2} 2^{-jp/2}  \; \Gamma \left(\frac{p}{2}\right) .\]
By Sterling's formula, $\Gamma (z) \ll \sqrt{\frac{2\pi}{z}} \left( \frac{z}{e}\right)^z$.
Thus, $\Gamma \left(\frac{p}{2}\right) \ll \sqrt{\frac{4\pi}{p}}\left(\frac{p}{2e}\right)^{\frac{p}{2}}$.
By arguing analogously for the imaginary parts, we obtain:
\[\mathbb{E}\left[ \max_{I \subseteq [k]} |S_I - \mathbb{E}[S_I]|^p\right] \ll C^{p} k^{\frac{p}{2}} p^{\frac{p}{2}} 2^{-jp/2},\]
where $C$ is a positive constant.
\end{proof}

Using the above lemma, we estimate $\mathbb{E}\left[ \left|\left| \{Z_n\}_{n=1}^{N_j}\right|\right|_{V^2}^2\right]$ as follows. We let $N'_j = 2^{m}$ be the smallest power of $2$ which is $\geq N_j$. We then decompose $[N'_j]$ into a family of dyadic intervals. More precisely, we define $\mathcal{F}$ to be the family of intervals of the form
\[ ((d-1)2^\ell, d2^\ell], \; \ell \in \{0, 1, \ldots, m\}, \; d \in \{1, \ldots, 2^{m-\ell}\}.\]
Now, for any interval $I'$, there are (at most) two intervals $I_l, I_r \in \mathcal{F}$ such that $I' \subseteq I_l \cup I_r$ and $|I_l \cup I_r| < 4|I'|$. Moreover, for any partition $\mathcal{P}$ of $[N_j]$, the number of times an $I \in \mathcal{F}$ is associated to an $I' \in \mathcal{P}$ is upper bounded by a constant. (This is as we have argued previously.)

We let $\Omega$ denote our probability space ($\omega \in \Omega$ corresponds to a specified value for each $Z_n$). Now, for a fixed $\omega \in \Omega$, we say an interval $I \subseteq \mathcal{F}$ is \emph{good} if:
\[ \max_{I' \subseteq I} |S_{I'} - \mathbb{E}[S_{I'}]|^2 \leq D 2^{-j} |I| \ln \ln (N),\]
where $D$ is a positive constant whose value we will specify later. Otherwise, we say $I$ is \emph{bad}. We let $\mathcal{P}$ denote the maximal partition (which depends on $\omega$). For each interval $I' \in \mathcal{P}$, we have (at most two) covering intervals $I_r, I_l \in \mathcal{F}$. We let $\mathcal{F}_{\mathcal{P}}$ denote the set of intervals in $\mathcal{F}$ which correspond to intervals in $\mathcal{P}$ (each $I \in \mathcal{F}$ corresponds to at most a constant number of intervals $I' \in \mathcal{P}$). We have:
\[\sum_{I' \in \mathcal{P}} \left| \sum_{n \in I} Z_n\right|^2 \ll \sum_{I \in \mathcal{F}_{\mathcal{P}}} \max_{I' \subseteq I} \left| \sum_{n \in I'} Z_n\right|^2.\]

We observe that
\[\sum_{\stackrel{I \in \mathcal{F}_{\mathcal{P}}}{I \text{ is good}}} \max_{I' \subseteq I} \left|\sum_{n \in I'} Z_n\right|^2 \ll \left|\sum_{c \in \mathcal{C}}c \right|^2 + D2^{-j} N_j \ln \ln (N) \ll \ln \ln (N) \sum_{c \in \mathcal{C}} |c|^2 + \left|\sum_{c \in \mathcal{C}}c \right|^2,\]
since each $|c|^2$ is between $2^{-j-1}$ and $2^{-j}$, and $|\mathcal{C}| = N_j$.
To see this, note that for each $I'$, $|S_{I'}|^2 \ll |S_{I'} - \mathbb{E}[S_{I'}]|^2 + |\mathbb{E}[S_{I'}]|^2$, and $|\mathbb{E}[S_{I'}]|^2 = \left| \frac{|I'|}{N_j} \sum_{c \in \mathcal{C}}c\right|^2$.

It only remains to bound the contribution of the intervals that are not good.
For this, we first prove the following lemma. For each interval $I \in \mathcal{F}$, we let $B(I)$ denote the event that $I$ is \emph{bad} (i.e. not good), and we let $1_{B(I)}$ denote its indicator function.
\begin{lemma}\label{lem:probbad}
For each $I \in \mathcal{F}$,
\[\mathbb{P}\left[ 1_{B(I)}\right] \ll \frac{1}{\ln (N)^4},\]
when $D$ is chosen to be a sufficiently large constant.
\end{lemma}

\begin{proof} By Chebyshev's inequality, for any $p>2$ we have
\begin{equation}\label{chebyshevresult}
\mathbb{P}\left[ 1_{B(I)}\right] = \mathbb{P}\left[ \max_{I' \subseteq I} |S_{I'} - \mathbb{E}[S_{I'}]|^2 > D 2^{-j} |I| \ln \ln (N)\right] \ll \frac{\mathbb{E}\left[ \max_{I'\subseteq I} |S_{I'} - \mathbb{E}[S_{I'}]|^p\right]}{\left(D 2^{-j} |I| \ln \ln (N)\right)^{p/2}}.
\end{equation}

We now rely on the following result of Ros\'{e}n \cite{Rosen}.
\begin{lemma}\label{lem:rosen}(Theorem 4 in \cite{Rosen}) Let $X_1, \ldots, X_k$ be samples drawn from a finite set of real numbers with replacement, and let $Z_1, \ldots, Z_k$ be samples drawn without replacement. Let $1 \leq n_1 < n_2 < \cdots < n_m$. For every convex, monotone function $\phi: \R \rightarrow \R$, we have
\[\mathbb{E}\left[ \max \left( \phi\left(\sum_{n=1}^{n_1} Z_n\right), \ldots, \phi\left(\sum_{n=1}^{n_m} Z_n\right)\right) \right] \leq \mathbb{E}\left[ \max \left( \phi\left(\sum_{n=1}^{n_1} X_n\right), \ldots, \phi\left(\sum_{n=1}^{n_m} X_n\right)\right) \right].\]
\end{lemma}
We want to apply this lemma to the function $f(x) := |x|^p$, but this is not monotone. Instead we define monotone, convex functions $f_1, f_2$ such that $|x|^p = f_1 (x) + f_2(x)$, namely setting $f_1(x) = (-x)^p$ for $x <0$ and equal to 0 otherwise, and $f_2(x) = x^p$ for $x >0$ and equal to 0 otherwise. We note that $|x|^p \geq f_1(x), f_2(x)$ always holds.

Without loss of generality, we consider $I$ equal to the interval of length $|I|$ starting at 1. Then, for some constant $H$, we have:
\[\mathbb{E}\left[ \max_{I'\subseteq I} |S_{I'} - \mathbb{E}[S_{I'}]|^p\right] \ll H^p \; \mathbb{E}\left[ \max_{1 \leq n \leq |I|} f_1\left(Re\left( S_n - \mathbb{E}[S_n] \right)\right)\right] + \]
\[\cdots + H^p \; \mathbb{E}\left[\max_{1 \leq n \leq |I|} f_2\left(Im\left(S_n - \mathbb{E}[S_n]\right)\right)\right].\]
Here, $S_n$ denotes the partial sum of $Z_1 + Z_2 + \cdots + Z_n$, $Re$ denotes the real part, $Im$ denotes the imaginary part, and there are four terms in this sum: one for each combination of $f_1,f_2$ and real and imaginary parts.

We can apply Lemma \ref{lem:rosen} to each of these four terms to replace the samples $Z_1, \ldots, Z_{|I|}$ taken without replacement with samples $X_1, \ldots, X_{|I|}$ taken with replacement. Now applying Lemma \ref{lem:complicated}, we have
\[\mathbb{P}\left[ 1_{B(I)}\right] \ll \frac{\tilde{H}^p |I|^{\frac{p}{2}} p^{\frac{p}{2}} 2^{-jp/2}}{ \sqrt{D}^p (\ln \ln (N))^{\frac{p}{2}} |I|^{\frac{p}{2}} 2^{-jp/2}} = \left(\frac{\tilde{H}}{\sqrt{D}}\right)^p p^{\frac{p}{2}} (\ln \ln (N))^{-\frac{p}{2}},\]
for some constant $\tilde{H}$.

Now, setting $p := \ln \ln (N)/e$, this is:
\[= \left( \frac{\tilde{H}}{\sqrt{D}}\right)^{\frac{\ln \ln (N)}{e}} \ln (N)^{-\frac{1}{2e}}.\]
We can then set $D$ large enough so that $\frac{\tilde{H}}{\sqrt{D}} < e^{-4e}$, and the lemma follows.
\end{proof}

We observe that the contribution of the bad intervals is upper bounded by
\begin{equation}\label{badsum}
\ll \sum_{ I \in \mathcal{F}} \mathbb{E}\left[ 1_{B(I)} \max_{I' \subseteq I} |S_{I'}|^2\right].
\end{equation}
We next apply H\"{o}lder's inequality with $q,r$ fixed to be constants such that $\frac{1}{r}+ \frac{1}{q} = 1$ and $\frac{4}{q} >2, r >1$. We then have that the above quantity is:
\[\ll \sum_{I \in \mathcal{F}} \left( \mathbb{E}[1_{B(I)}]\right)^{\frac{1}{q}} \left( \mathbb{E}\left[ \max_{I' \subseteq I} |S_{I'}|^{2r}\right]\right)^{\frac{1}{r}}.\]

By Lemma \ref{lem:probbad}, we know that
\[\left( \mathbb{E}[1_{B(I)}]\right)^{\frac{1}{q}} \ll (\ln (N))^{-2}.\]
We also know that for each $I'$, $|\mathbb{E}[S_{I'}]|^2 \ll \left(\frac{|I'|}{N_j}\right)^2 \left|\sum_{c \in \mathcal{C}} c\right|^2 \ll \frac{|I'|}{N_j} \left|\sum_{c \in \mathcal{C}} c\right|^2$. When we sum these up over all $\mathcal{I} \in \mathcal{F}$, we obtain $\ll \ln (N) \left| \sum_{c\in \mathcal{C}} c \right|^2$. Now multiplying by $ \ln(N)^{-2}$, we obtain a contribution which is $o\left( \left|\sum_{c\in \mathcal{C}} c\right|^2\right)$. Thus, it only remains to bound
\[(\ln (N))^{-2}\sum_{I \in \mathcal{F}} \left( \mathbb{E}\left[ \max_{I' \subseteq I} |S_{I'}-\mathbb{E}[S_{I'}]|^{2r}\right]\right)^{\frac{1}{r}}.\]

Similarly to our above arguments, we define convex, monotone functions $f_1, f_2: \R \rightarrow \R$ such that $f_1(x) + f_2(x) = |x|^{2r}$. More precisely, we set $f_1(x) = (-x)^{2r}$ when $x <0$ and equal to 0 otherwise, while we set $f_2(x) = x^{2r}$ when $x >0$ and equal to 0 otherwise. Now, again applying Lemma \ref{lem:rosen}, it suffices to bound e.g.
\[ \sum_{I \in \mathcal{F}} \left(\mathbb{E}\left[\max_{1 \leq n \leq |I|} f_1(Re(S_n - \mathbb{E}[S_n]))\right]\right)^{\frac{1}{r}},\]
where $S_n$ is now the partial sum $X_1 + \cdots + X_n$, where each $X_k$ is a sample from $\mathcal{C}$ taken \emph{with} replacement. (We must also bound the analogous quantities for other combinations of $f_1, f_2$ and $Re, Im$, but these will follow via the same argument.)

We now apply Lemma \ref{lem:Doob} to obtain that the above quantity is
\[ \ll \sum_{I \in \mathcal{F}} \left(\mathbb{E}\left[ \max_{1 \leq n \leq |I|}|Re(S_n-\mathbb{E}[S_n])|^{2r}\right]\right)^{\frac{1}{r}} \ll \sum_{I \in \mathcal{F}} \left( \mathbb{E}[|Re(S_{I}- \mathbb{E}[S_I])|^{2r}]\right)^{\frac{1}{r}}.\]
Next applying Lemma \ref{lem:Rosenthal}, we see that this is
\[\ll \sum_{I \in \mathcal{F}} \max \left\{ \left( \sum_{n=1}^{|I|} \mathbb{E}[|\tilde{X}_n|^{2r}]\right)^{\frac{1}{r}}, \sum_{n=1}^{|I|} \mathbb{E}[|\tilde{X}_n|^2]\right\},\]
where $\tilde{X}_n$ is defined to be an (independent, uniform) sample from $\mathcal{C}$ with replacement, recentered to be mean zero. In other words, $\tilde{X}_n = X_n - \mathbb{E}{X_n}$.
Now, since $r >1$, both of the quantities in this maximum are $\ll |I| 2^{-j}$. Hence, we have:
\[\ll \sum_{ I \in \mathcal{F}} |I| 2^{-j} \ll \ln (N) \sum_{c \in \mathcal{C}} |c|^2.\]
Multiplying this by our bound $(\ln (N))^{-2}$ for the probability of each $I$ being bad, we see that this is $o\left( \sum_{c \in \mathcal{C}} |c|^2\right)$. This completes the proof of Lemma \ref{lem:shortop}.

\end{proof}

Combining Lemma \ref{lem:shortop} with (\ref{blockgoal}), we obtain Theorem \ref{mod1Perm}.

\end{proof}

\section{Refinements of Theorem \ref{main} for Certain Structured ONS}\label{sec7}

In this section, we briefly outline how Theorem \ref{main} can be improved for more restrictive classes of ONS, using the methods employed in proving Theorem \ref{varRad}. We consider an ONS such that for $f$ in the span of the system, we have $||f||_{L^p} \leq C_p ||f||_{L^2}$ for $p>2$, where $C_p$ is a constant depending only on $p$. Such systems arise naturally, for example, as the restriction of the trigonometric system to certain arithmetic subsets ($\Lambda(p)$ sets). We will use the fact that a maximal form of this hypothesis can be obtained from a very general theorem of Christ and Kiselev \cite{ChristKiselev}.

\begin{theorem} Let $\{ \phi_n \}_{n=1}^{\infty}$ be an ONS such that for $f$ in the span of the system, we have $||f||_{L^p} \leq C_p ||f||_{L^2}$ for some $p>2$. Then
\begin{equation}
||\mathcal{M}f||_{L^p} \ll_{\delta} C_{p} ||f||_{L^2}
\end{equation}
as long as $p>\delta>2$.
\end{theorem}

This last condition implies that the implicit constant is uniform for large $p$. Using this and the arguments in the proof of Theorem \ref{varRad}, one can obtain the following:

\begin{theorem}
Let $\{ \phi_n \}_{n=1}^{\infty}$ be a ONS such that if $f$ is in the span of the system, then $||f||_{L^p} \ll C_{p} ||f||_{L^2}$ for some $p>2$. We then have that
\[||f||_{L^2(V^2)} \ll_{p} \ln^{1/p}(|A|)||f||_{L^2}. \]
where the coefficients of $f$ are supported a finite index set $A$.
\end{theorem}

We briefly sketch the proof. We note that if $||\mathcal{M}f||_{L^2} \ll ||f||_{L^2}$ holds, then this theorem follows for $p=2$. However, this is in general not true and by the sharpness of Theorem \ref{main}, the best one can hope for in the general case is a factor of $\ln(|A|)$ in place of $\ln^{1/2}(|A|)$. The proof follows the same setup as the proof of Theorem \ref{varRad}. We define a bad event for some interval $J$ to be the event that $ |\tilde{S}_J| \gg \ln^{1/p}(|A|) (M(J))^{1/2}$ (here $M(J)$ is defined to be the sum of $a_n^2$ over $n \in J$, where the $a_n$'s are the coefficients of $\phi_n$ in the expansion of $f$).  It is easy to see that the contribution from the good events are of an acceptable order and it suffices to bound the bad events. The argument is essentially the same as the proof of Theorem \ref{varRad}, with the exception that we use the following estimate:
\[\int_{\T} |1_{B(\tilde{J})} \tilde{S}_{\tilde{J}}|^2 \leq \left(\int_{\T} 1_{B(\tilde{J})}  \right)^{1/(p/2)'} \left( \int_{\T} |\tilde{S}_{\tilde{J}}|^p \right)^{(2/p)}.  \]
(Here, $(p/2)'$ denotes the conjugate exponent of $p/2$.)

We now estimate $\int_{\T}|\tilde{S}_{\tilde{J}}|^p  \ll C_{p}^p \; \left(\int_{\T} |S_{\tilde{J}}|^2 \right)^{p/2} \ll C_{p}^p \; (M(\tilde{J}))^{p/2}$. Hence $ \left(\int_{\T} |\tilde{S}_{\tilde{J}}|^p \right)^{(2/p)} \ll C_{p}^2 M(\tilde{J})$.
Next, by Chebyshev's inequality, \[\int_{\T}  1_{B(\tilde{J})} \leq \frac{\int_{\T} |\tilde{S}_{\tilde{J}}|^p }{ \left(  \ln^{1/p} (|A|) (M(\tilde{J}))^{1/2} \right)^{p}  } \leq \frac{ C_{p}^p}{\ln(|A|)}.\]

Hence (using $1/(p/2)'= \frac{p-2}{p}$), we have $\left(\int_{\T}  1_{B(\tilde{J})} \right)^{(p-2)/p} \ll \frac{ C_{p}^{p-2} }{\ln^{(p-2)/p}(|A|)}$. This yields
\[\int_{\T} |1_{B(J)} \tilde{S}_{\tilde{J}}|^2  \leq \left(\int_{\T} 1_{B(\tilde{J})}  \right)^{1/(p/2)'} \left(\int_{\T} |\tilde{S}_{\tilde{J}}|^p \right)^{(2/p)} \ll \frac{C_{p}^{p} M(\tilde{J})}{\ln^{(p-2)/p} (|A|)}.  \]

Now we sum this quantity over $\ln(|A|)$ levels, each with the sum of $M(\tilde{J})$ summing to $1$. Hence the contribution from the bad events to the quantity we wish to estimate is $O( \ln^{2/p}(|A|))$. This is exactly the order we wish to show.

Finally, we observe that:
\begin{theorem}Let $\{ \phi_n \}_{n=1}^{\infty}$ be an ONS such that if $f$ is in the span of the system, then $||f||_{L^p} \ll \sqrt{p}||f||_{L^2}$ (for all $p>2$). Then
\[||f||_{L^2(V^2)} \ll \sqrt{\ln\ln(|A|)}||f||_{L^2}, \]
where the coefficients of $f$ are supported on the index set $A$.
\end{theorem}

This is proved using the same arguments sketched for the previous theorem, however now we have freedom to optimize over the choice of $p$ we use. The optimum occurs with a choice of $p$ about $c e^{-1}\ln\ln(N)$. Essentially the same argument is given in detail in the proof of Theorem \ref{mod1Perm} for random permutations (see the proof of Lemma \ref{lem:probbad}). Here it is important that the constants in the Christ-Kiselev theorem are uniformly bounded for large $p$.

The above theorem can be applied to systems formed by Sidon subsets of the trigonometric system, since the hypothesis of this theorem characterizes Sidon sets (when applied to subsets of the trigonometric system) by a theorem of Pisier \cite{Pisier} (see also \cite{Rudin}).

\section{Variational Estimates for the $V^p$ Operator}

\subsection{Notation}

Let $\Gamma: \R \rightarrow \R^{+}$ be a convex symmetric function, increasing on $\R^{+}$ and tending to infinity at infinity such that $\Gamma(0)=0$. Then the Orlicz space norm associated to $\Gamma$ is defined as

\[||f||_{\Gamma}:= \min \left\{\lambda : \int_{\T}  \Gamma \left(\frac{f(x)}{\lambda} \right) dx \leq 1 \right\}.\]

The fact that this norm satisfies the triangle inequality is an easy exercise using Jensen's inequality. We refer the reader to \cite{KR} for the general theory of these spaces. Following \cite{Bour}, we will be interested in $\Gamma := \Gamma_{K}$ defined as follows

\[\Gamma_{K}(t):= \left\{
                   \begin{array}{ll}
                     |t|^{5/2}, & \hbox{$|t|\leq K$} \\
                     \frac{5}{4}K^{1/2}t^2 - \frac{1}{4}K^{5/2}, & \hbox{$|t| \geq K$}
                   \end{array}
                 \right. .
\]
Later we will also use
\[\gamma_{K}(t) :=\left\{
                    \begin{array}{ll}
                    |t|^{1/2}, & \hbox{$|t|\leq K$} \\
                    K^{1/2}, & \hbox{$|t| \geq K$}
                \end{array}
                \right. .
\]

We note that $t^2 \gamma_K(t) \leq \Gamma_K(t)$ for all $t$. We state some other basic properties that we will need.

\begin{lemma}\label{lem:pconvex} Let $2= p$. Then $||\cdot||_{\Gamma_K}$ is $p$-convex. That is, for any functions $f_1, \ldots, f_k$ from $\T$ to $\R$,
\[\left|\left| \left(\sum_{i=1}^k |f_{i}|^p\right)^{1/p} \right|\right|_{\Gamma_K} \leq  \left(\sum_{i=1}^k  || f_i ||_{\Gamma_{K}}^{p}\right)^{1/p}.\]

\end{lemma}

\begin{proof}

Let $\Gamma_{K,1/p}(t):= \Gamma_{K}(t^{1/p})$, which we observe is still convex (we have used that $p= 2 $ here). Since $\Gamma_{K,1/p}(t)$ is convex, we can use it to form an Orlicz space norm.
We observe that
\[\left|\left| \left(\sum_{i=1}^k |f_{i}|^p\right)^{1/p} \right|\right|_{\Gamma_K} = \min \left\{\lambda : \int_{\T}  \Gamma_{K} \left(\frac{\left(\sum_{i=1}^k |f_{i}(x)|^p\right)^{1/p} }{\lambda} \right) dx \leq 1 \right\}\]
\[ = \min \left\{\lambda : \int_{\T}  \Gamma_{K,1/p} \left(\frac{ \sum_{i=1}^k |f_{i}(x)|^p }{\lambda^p} \right) dx \leq 1 \right\} = \left|\left| \sum_{i=1}^k |f_{i}|^p \right|\right|_{\Gamma_{K,1/p}}^{1/p}  \]
\[\leq \left(\sum_{i=1}^k \left|\left| |f_{i}|^p \right|\right|_{\Gamma_{K,1/p}} \right)^{1/p} = \left(\sum_{i=1}^k  || f_i ||_{\Gamma_{K}}^{p}\right)^{1/p}.\]
The inequality here follows from the triangle inequality for $||\cdot||_{\Gamma_{K,1/p}}$.
\end{proof}

\subsection{Proof of Theorem \ref{varVp}}
We now prove:
\begin{thm10} Let $p>2$ and $\{\phi_{n}\}_{n=1}^{N}$ be an orthonormal system such that $||\phi_{n}||_{L^\infty}\leq C$ for all $n$. There exists a permutation $\pi:[N]\rightarrow [N]$ such that the orthonormal system
$\{\psi_n := \phi_{\pi(n)}\}_{n=1}^{N}$ satisfies
\begin{equation}\label{varperm}
||f||_{L^{2}(V^{p})} \ll_{C,p} \ln\ln(N)||f||_{L^2}
\end{equation}
for all $f= \sum_{n=1}^{N} a_n\psi_n(x)$.
\end{thm10}

Our starting point is the inequality (3.21) of \cite{Bour}:

\begin{theorem}\label{bourgain} Let $\{\phi_n\}_{n=1}^N$ be an orthonormal system with $||\phi_n||_{L^\infty} \leq C$ for all $n$.
Then there exists a permutation $\pi: [N] \rightarrow [N]$ such that for all subintervals $I$ of $[N]$ and all real values $a_1, \ldots, a_N$, the orthonormal system $\{\psi_n := \phi_{\pi(n)}\}_{n=1}^N$ satisfies:

\begin{equation}\label{permIneq}
\left|\left| \sum_{n \in I} a_{n}\psi_{n} \right|\right|_{\Gamma_{N/|I|}} \ll_C \ln^{3/4}(N) \left( \sum_{n \in [I]} a_n^2\right)^{1/2}.
\end{equation}

\end{theorem}

We will need a variational form of this inequality. This is easily achieved using a Rademacher-Menshov argument.

\begin{lemma}\label{lem:perm} With the notation as above, we have that
\begin{equation}\label{permIneq}
 \bigg| \bigg|   ||  \{a_{n}\psi_{n}\}_{n \in I} ||_{V^2}   \bigg|\bigg|_{\Gamma_{N/|I|}}\ll_C \ln^{7/4}(N) \left( \sum_{n \in I} a_n^2\right)^{1/2}
\end{equation}
for all $I \subseteq [N]$ and all real sequences $a_1, \ldots, a_N$.
\end{lemma}

\begin{proof}
As in section \ref{sec:rad}, we assume (without loss of generality) that $I = [2^{\ell}]$ for some $\ell$ and we define the intervals $I_{k,i} := (k2^{i}, (k+1)2^i]$ for $0 \leq i \leq \ell$ and $0 \leq k \leq 2^{\ell-i}-1$. For each $J \subseteq I$, we can express $J$ as a disjoint union of intervals $I_{k,i}$, where the union contains at most two intervals of each size. As in (\ref{pointwise}), we then observe for each $x \in \T$:
\[|| \{a_{n}\psi_{n}\}_{n \in I} ||_{V^2}(x) \ll \sum_{i=0}^\ell \sqrt{\sum_{k=0}^{2^{\ell-i}-1} \left( \sum_{n \in I_{k,i}} a_n \psi_n(x)\right)^2}. \]

By the triangle inequality for the Orlicz norm, we then have
\[ \bigg| \bigg|   ||  \{a_{n}\psi_{n}\}_{n \in I} ||_{V^2}   \bigg|\bigg|_{\Gamma_{N/|I|}} \ll \sum_{i=0}^\ell \left|\left| \sqrt{\sum_{k=0}^{2^{\ell-i}-1} \left( \sum_{n \in I_{k,i}} a_n \psi_n(x)\right)^2}\right|\right|_{\Gamma_{N/|I|}}.\]

Applying Lemma \ref{lem:pconvex}, this is
\[\leq \sum_{i=0}^\ell \sqrt{ \sum_{k=0}^{2^{\ell-i}-1} \left|\left|\sum_{n \in I_{k,i}} a_n \psi_n(x) \right|\right|^2_{\Gamma_{N/|I|}}}.\]
By Theorem \ref{bourgain}, we obtain
\[ \ll_C \ln^{3/4} (N)\sum_{i=0}^\ell \sqrt{ \sum_{k=0}^{2^{\ell-i}-1} \sum_{n \in I_{k,i}} a_n^2}  = \ln^{3/4} (N) \sum_{i=0}^\ell \sqrt{\sum_{n \in I} a_n^2} = \ln^{7/4} (N) \sqrt{\sum_{n \in I} a_n^2}.\]

\end{proof}

We now prove Theorem \ref{varVp}. We assume (without loss of generality) that $\sum_{n=1}^N a_n^2 =1$. As in Section \ref{sec:probability}, we consider decomposing $[N]$ into a family of subintervals according to mass, defined with respect to the $a_n$'s. We recall that the mass of an arbitrary subinterval $I$ is defined to be $M(I):= \sum_{n \in I} a_n^2$. We define the intervals $I_{k,s}$ for $1 \leq s \leq 2^k$ and points $i_{k,s}$ as in Section \ref{sec:probability}. We refer to the intervals $I_{k,s}$ for $1 \leq s \leq 2^k$ as the admissible intervals on level $k$, and the points $i_{k,s}$ (as $s$ ranges) as the admissible points on level $k$.
We note that any interval $I \subseteq [N]$ can be expressed as a union of intervals of the form $I_{k,s}$ and points $i_{k,s}$, where there are at most two intervals and two points for each value of $k$ (this follows analogously to the proof of Lemma \ref{lem:binarydecomp}). This decomposition is obtained by first taking the intervals $I_{k,s}$ and points $i_{k,s}$ contained in $I$ with the smallest value of $k$. (There are at most 2 of each, otherwise $I$ would contain an admissible interval or point for a smaller $k$ value.) These ``components" of $I$ on level $k$ form an interval, and when we remove this from $I$, we are left with a left part and a right part. Each part can then be decomposed as union of intervals $I_{k,s}$ and points $i_{k,s}$ for higher values of $k$, and each of the two unions contains at most one interval and one point on each level.

We let $\pi:[N] \rightarrow [N]$ be the permutation as in Lemma \ref{lem:perm}, and $\psi_n := \phi_{\pi(n)}$.
We fix an $x \in \T$. The value of
\[\left| \left| \{a_n \psi_n(x)\}_{n=1}^N \right| \right|_{V^p}\]
is achieved by some partition $\mathcal{P}$ of $[N]$. Each $I \in \mathcal{P}$ can be expressed as a union of intervals of the form $I_{k,s}$ and points $i_{k,s}$, and we denote the set of these intervals and points by $T_I$ and $t_I$ respectively. We recall that each of $T_I$ and $t_I$ will have at most two intervals or points (respectively) on each level. We also note that each admissible interval will appear in this union for at most one $I \in \mathcal{P}$.

We fix a positive constant $c$ (depending on $p$) such that $c > \max \{\frac{35}{4}\left(\frac{1}{2} - \frac{1}{p}\right)^{-1}, 9\}$ (this is possible because $p >2$). We define $k^* :=c  \ln \ln (N)$ (more precisely, $k^*$ is the nearest integer greater than $c \ln \ln (N)$). Now, for each $I \in \mathcal{P}$, all of the intervals in $T_I$ and points in $t_I$ on levels greater than $k^*$ are contained in the two intervals $I_{k^*,s_\ell}$ and $I_{k^*,s_r}$ on level $k^*$, where $s_\ell$ is one less than the $s$ value for the leftmost interval $I_{k^*,s}$ in $T_I$, and $s_r$ is one more than the $s$ value for the rightmost interval $I_{k^*,s}$ in $T_I$. We will use $k^*$ as a cutoff threshold: we handle the intervals and points at levels $\leq k^*$ directly and handle the intervals and points at levels $> k^*$ using the fact that they are contained in $I_{k^*,s_\ell}, I_{k^*,s_r}$. We define $T'_I$ to be the subset of intervals in $T_I$ on levels $\leq k^*$ and $t'_I$ to be the subset of points in $t_I$ on levels $\leq k^*$.

Now, $\left| \left| \{a_n \psi_n(x)\}_{n=1}^N \right| \right|_{V^p}$ is equal to:
\[ \left(\sum_{I \in \mathcal{P}} \left( \sum_{n \in I} a_n \psi_n(x)\right)^p\right)^{1/p} =\]\[
\left(\sum_{I \in \mathcal{P}} \left( \sum_{J \in T'_I} \sum_{n \in J} a_n \psi_n(x) + \sum_{J \in T_I\backslash T'_I} \sum_{n \in J} a_n \psi_n(x)
+\sum_{n \in t'_I} a_n \psi_n(x) + \sum_{n \in t_I \backslash t'_I} a_n \psi_n(x)\right)^p\right)^{1/p}.\]
Applying the triangle inequality for the $\ell_p$-norm, this is:
\begin{eqnarray}\label{decomp1}
\nonumber
  &\leq& \left( \sum_{I \in \mathcal{P}} \left( \sum_{J \in T'_I} \sum_{n \in J} a_n \psi_n(x)\right)^p\right)^{1/p} + \left(\sum_{I \in \mathcal{P}} \left( \sum_{n \in t'_I} a_n \psi_n(x)\right)^p\right)^{1/p} \\
  &+ & \left(\sum_{I \in \mathcal{P}} \left( \sum_{J \in T_I\backslash T'_I} \sum_{n \in J} a_n \psi_n(x) + \sum_{n \in t_I\backslash t'_I} a_n \psi_n(x)\right)^p \right)^{1/p}
\end{eqnarray}


We consider the second of these three terms. Since $p \geq 2$, we have
\[\left(\sum_{I \in \mathcal{P}} \left( \sum_{n \in t'_I} a_n \psi_n(x)\right)^p\right)^{1/p} \leq \left(\sum_{I \in \mathcal{P}} \left(\sum_{n \in t'_I} a_n \psi_n(x)\right)^2\right)^{1/2}.\]
For each $k \leq k^*$, we let $\ell_k$ denote the set of admissible points on level $k$. Since each $t'_I$ contains at most 2 points in each $\ell_k$, we can apply the triangle inequality to obtain
\[\left(\sum_{I \in \mathcal{P}} \left(\sum_{n \in t'_I} a_n \psi_n(x)\right)^2\right)^{1/2} \ll \sum_{k=0}^{k^*}\left(\sum_{n \in \ell_k} (a_n\psi_n(x))^2\right)^{1/2}.\]

Now, by the triangle inequality for the $L^2$ norm and the fact that $\int_{\T} a_n^2 \psi_n^2(x) dx= a_n^2$ for all $n$, we have
\[\left| \left| \sum_{k=0}^{k^*} \left(\sum_{n \in \ell_k} (a_n\psi_n(x))^2\right)^{1/2}\right|\right|_{L^2} \ll_p \ln \ln (N).\]
To see this, recall that $\sum_{n=1}^N a_n^2 = 1$, so $\sum_{n \in \ell_k} a_n^2 \leq 1$ for each $k$, and $k^* \ll_p \ln \ln (N)$.

It remains to bound the first and third terms in (\ref{decomp1}). We consider the first term. For each $k$, we let $\mathcal{L}_k$ denote the set of admissible intervals $I_{k,s}$ as $s$ ranges from 1 to $2^k$ (i.e. the admissible intervals on level $k$).
Then, by triangle inequality for the $\ell^2$ norm and the fact that $p \geq 2$,
\[\left( \sum_{I \in \mathcal{P}} \left( \sum_{J \in T'_I} \sum_{n \in J} a_n \psi_n(x)\right)^p\right)^{1/p} \leq \left( \sum_{I \in \mathcal{P}} \left( \sum_{J \in T'_I} \sum_{n \in J} a_n \psi_n(x)\right)^2\right)^{1/2}\]
\[\leq \sum_{k=0}^{k^*} \left(\sum_{I \in \mathcal{P}} \left (\sum_{J \in T'_I \cap \mathcal{L}_k} \sum_{n \in J} a_n\psi_n(x)\right)^2\right)^{1/2}\]
\[\leq \sum_{k=0}^{k^*} \left( \sum_{J \in \mathcal{L}_k} \left( \sum_{n\in J} a_n \psi_n(x)\right)^2\right)^{1/2}.\]

Now, using the triangle inequality for the $||\cdot ||_{L^2}$ norm, we have:
\[\left|\left| \sum_{k=0}^{k^*} \left( \sum_{J \in \mathcal{L}_k} \left( \sum_{n\in J} a_n \psi_n(x)\right)^2\right)^{1/2} \right| \right|_{L^2} \leq
\sum_{k=0}^{k^*} \left|\left| \left( \sum_{J \in \mathcal{L}_k} \left( \sum_{n\in J} a_n \psi_n(x)\right)^2\right)^{1/2} \right|\right|_{L^2}\]
\[ = \sum_{k=0}^{k^*} \left( \sum_{J \in \mathcal{L}_k} \int_{\T}\left( \sum_{n \in J} a_n \psi_n(x)\right)^2 dx \right)^{1/2}\]
\[ = \sum_{k=0}^{k^*} \left( \sum_{J \in \mathcal{L}_k} M(J)\right)^{1/2}  \ll_p \ln \ln (N),\]
since $\sum_{J \in \mathcal{L}_k} M(J) =1$ for each $k$, and $k^* \ll_p \ln \ln (N)$.

We are thus left with the third term of (\ref{decomp1}). For each $I \in \mathcal{P}$, we consider the union of the intervals and points in $T_I \backslash T'_I$ and $t_I \backslash t'_I$. This can alternatively be described as a union of at most two intervals $J_\ell$ and $J_r$, where each of $J_\ell, J_r$ is a subinterval of $I_{k^*,s}$ for some $s$. To see this, recall that $I$ is decomposed into a union of admissible intervals and points by taking the admissible intervals and points contained in $I$ for the earliest level where this set is non-empty. The remaining left and right parts of $I$ are then decomposed separately. If the minimal $k$ is $\leq k^*$, then $J_\ell$ is the union of the intervals/points in the decomposition of the left part that fall beyond level $k^*$, and $J_r$ is the same for the right part. If the minimal $k$ is $> k^*$, then in fact all of $I$ is contained in some admissible interval on level $k^*$, and we can take $J_\ell$ to be this interval and $J_r$ to be empty. We then rewrite the quantity we wish to bound as:
\[\left( \sum_{I \in \mathcal{P}} \left( \sum_{n \in J_\ell} a_n \psi_n(x) + \sum_{n \in J_r} a_n\psi_n(x)\right)^p \right)^{1/p}.\]
Applying the simple fact that $(a+b)^p \leq 2^p(a^p + b^p)$ for all non-negative real numbers $a$ and $b$, we see this is
\[\ll \left( \sum_{I \in \mathcal{P}} \left( \sum_{n \in J_\ell} a_n \psi_n(x)\right)^p + \left(\sum_{n \in J_r} a_n \psi_n(x)\right)^p \right)^{1/p}.\]

Now we observe that we are summing the values $a_n \psi_n(x)$ over disjoint intervals, each of which is contained in $I_{k^*, s}$ for some $s$. Thus, this quantity is upper bounded by:
\[\leq \left( \sum_{1 \leq s \leq 2^{k^*}} \left| \left| \{a_n \psi_n(x)\}_{n \in I_{k^*,s}} \right|\right|^p_{V^p}\right)^{1/p}.\]
Therefore, it suffices to bound
\[\left|\left| \left( \sum_{1 \leq s \leq 2^{k^*}} \left| \left| \{a_n \psi_n(x)\}_{n \in I_{k^*,s}} \right|\right|^p_{V^p}\right)^{1/p}\right|\right|_{L^2}.\]

%
%

For each $s$ from 1 to $2^{k^*}$, we define disjoint sets $G_s, B_s$ such that $G_s \cup B_s = \mathbb{T}$.  We define $G_s$ to be $x \in \mathbb{T}$ such that  $||\{a_n\psi_n(x)\}_{n \in I_{k^*,s}}||_{V^p} \leq 2^{- c\ln\ln(N) /p}$ and $B_s$ to be the complement. By two applications of the triangle inequality (first in the $\ell^p$ norm and then in the $L^2$ norm), we have
\[
  \left| \left|    \left( \sum_{s=1}^{2^{k^*}} ||\{a_n \psi_n(x)\}_{n \in I_{k^*,s}} ||_{V^{p}}^p \right)^{1/p} \right|\right|_{L^2} \ll   \left| \left|    \left( \sum_{s=1}^{2^{k^*}} 1_{G_s} ||\{a_n\psi_n(x)\}_{n \in I_{k^*,s}} ||_{V^{p}}^p \right)^{1/p} \right|\right|_{L^2} \]\[+  \left| \left|    \left( \sum_{s=1}^{2^{k^*}} 1_{B_s} ||\{a_n \psi_n(x)\}_{n \in I_{k^*,s}}||_{V^{p}}^p \right)^{1/p} \right|\right|_{L^2}.
\]

Using that $||\{a_n \psi_n(x)\}_{n \in I_{k^*,s}}||_{V^p}^{p} \ll 2^{-c\ln\ln(N)}$ for $x \in G_s$, we have that the first term is $O(1)$ (from the fact that there are at most $2^{c\ln \ln(N)}$ terms in the sum). We now estimate

\[\left| \left|    \left( \sum_{s=1}^{2^{k^*}} 1_{B_s}(x) ||\{a_n \psi_n(x)\}_{n \in I_{k^*,s}}||_{V^{p}}^p \right)^{1/p} \right|\right|_{L^2} \ll \left| \left|    \left( \sum_{s=1}^{2^{k^*}} 1_{\tilde{B}_s}(x) ||\{a_n \psi_n(x)\}_{n \in I_{k^*,s}} ||_{V^{2}}^2 \right)^{1/2} \right|\right|_{L^2} \]
\begin{equation}\label{V2lastterm}
 \ll \left(    \sum_{s=1}^{2^{k^*}} ||1_{\tilde{B}_s}(x) || \{a_n \psi_n(x)\}_{n\in I_{k^*,s}}||_{V^2} ||_{L^2}^2 \right)^{1/2},
 \end{equation}
where $\tilde{B}_s$ is the set of $x \in \T$ such that $||\{a_n \psi_n(x)\}_{n \in I_{k^*,s}}||_{V^2} \geq 2^{-c\ln\ln(N) /p}$, and we have used the fact that $B_s \subseteq \tilde{B}_s$.

We now consider two cases. First, we consider the set $S_{\text{big}}$ of $s$ values where $|I_{k^*,s}| \geq N2^{-7\ln \ln (N)}$. Clearly, there can be at most $2^{7 \ln \ln (N)}$ such intervals. Now we bound the contribution to (\ref{V2lastterm}) above from these big intervals as
\[\left(    \sum_{s \in S_{\text{big}}}^{2^{k^*}} ||1_{\tilde{B}_s}(x) || \{a_n \psi_n(x)\}_{n\in I_{k^*,s}}||_{V^2} ||_{L^2}^2 \right)^{1/2} \ll
\left(    \sum_{s \in S_{\text{big}}}^{2^{k^*}} || \{a_n \psi_n(x)\}_{n\in I_{k^*,s}}||_{L^2(V^2)}^2 \right)^{1/2}.\]
Recalling that $|| \{a_n \psi_n(x)\}_{n\in I_{k^*,s}}||_{L^2(V^2)}^2 \ll \ln^{2}(N)2^{-c \ln\ln(N)}$ (from Lemma \ref{varRM1}, since $M(I_{k^*,s}) \leq 2^{-k^*}$ for all $s$) and that there are at most $2^{7 \ln \ln (N)}$ values of $s \in S_{\text{big}}$, we have that the above is
\[ \ll \left( 2^{7\ln \ln (N)}  \ln^{2}(N)2^{-c \ln\ln(N)} \right)^{1/2} \ll 1.\]
Here we have used that $9 \leq c$.
It now suffices to consider the values of $s$ such that $|I_{k^*,s}| \leq N 2^{-7\ln \ln (N)}$.

We define $\gamma_{*}=\gamma_{2^{7\ln\ln(N)}}$.
For any real numbers $\epsilon >0$, $\lambda>1$, and $a \geq \epsilon$, we have $\frac{\gamma_{*}(\lambda^{-1}a) }{\gamma_{*}(\lambda^{-1}\epsilon) } \geq 1$. We set $\epsilon := 2^{-c\ln \ln (N)/p}$. Now, for all $x \in \tilde{B}_s$, we have:
\begin{equation}\label{easy}
\left| \left|   \{a_n\psi_n(x)\}_{n \in I_{k^*,s}} \right|\right|_{V^2}^2 \leq \left| \left|   \{a_n\psi_n (x)\}_{n \in I_{k^*,s}} \right|\right|_{V^2}^2 \frac{\gamma_{*}(\lambda^{-1} ||\{a_n\psi_n(x)\}_{n \in I_{k^*,s}}||_{V^2}) }{\gamma_{*}(\lambda^{-1}\epsilon) } .
\end{equation}

We recall that $M(I_{k^*,s}) \leq 2^{-c\ln\ln(N)}$ for each $s$. Analogously to $\gamma_{*}$, we define $\Gamma_{*} := \Gamma_{2^{7\ln \ln (N)}}$.
Now, for any $\lambda > 1$:
\[ \int_{\tilde{B}_{s}}\left| \left|   \{a_n \psi_n(x)\}_{n \in I_{k^*,s}} \right|\right|_{V^p}^2 dx \leq  \lambda^2 \int_{\tilde{B}_{s}} \gamma_{*}\left(\frac{\epsilon}{\lambda} \right)^{-1} \Gamma_{*}(\lambda^{-1}||\{a_n \psi_n(x)\}_{n \in I_{k^*,s}}||_{V^2}) dx. \]
This follows from (\ref{easy}) and the definitions of $\gamma_{*}$ and $\Gamma_{*}$ (recall also that $t^2 \gamma_{*}(t) \leq \Gamma_{*}(t)$ for all $t$).

Since $\frac{N}{|I_{k^*,s}|} \geq 2^{7\ln \ln (N)}$ and the value of $||\cdot ||_{\Gamma_K}$ increases as $K$ increases, we can apply Lemma \ref{lem:perm} to obtain
 \[ \left| \left|   | |  \{a_{n}\psi_{n}\}_{n \in I_{k^*,s}} ||_{V^2}   \right|\right|_{\Gamma_{*}}\leq D \ln^{7/4}(N) \left( \sum_{n \in I_{k^*,s}} a_n^2\right)^{1/2}\]
for all $s$ such that $|I_{k^*,s}| \leq N2^{-\frac{7}{2}\ln\ln (N)}$, where $D$ is some fixed constant (depending on $C$).

We see that for $\lambda := D \ln^{7/4}(N) 2^{- \frac{c \ln\ln(N)}{2}}$, we have $\int_{\mathbb{T}} \Gamma_{*}(\lambda^{-1}|| \{a_n \psi_n(x)\}_{n \in I_{k^*,s}} ||_{V^2}) dx \ll 1$. Therefore:
\begin{equation}\label{mark}
\int_{\tilde{B}_{s}}\left| \left| \{a_n \psi_n(x)\}_{n \in I_{k^*,s}} \right|\right|_{V^p}^2 dx \ll  \ln^{7/2}(N) 2^{- c \ln\ln(N)} \gamma_{*}\left(\frac{\epsilon}{\lambda} \right)^{-1}.
\end{equation}
We consider the quantity $\gamma_{*}\left(\frac{\epsilon}{\lambda}\right)^{-1}$. We observe:
\begin{equation}\label{quantity}
\frac{\epsilon}{\lambda} = (D^{-1})2^{\ln \ln (N) \left(-c/p +c/2-7/4\right)}.
\end{equation}
Now, if (\ref{quantity}) is $\geq 2^{7 \ln \ln (N)}$, we will have
\[\gamma_{*}\left(\frac{\epsilon}{\lambda}\right)^{-1} = 2^{-7/2 \ln \ln (N)}.\]
If (\ref{quantity}) is $< 2^{7 \ln \ln (N)}$, we will have
\[\gamma_{*} \left(\frac{\epsilon}{\lambda}\right)^{-1} = D^{1/2}2^{\ln \ln (N)(7/8-c/4+c/2p)}.\]
We note that $\frac{7}{8} - \frac{c}{4} + \frac{c}{2p} \leq -\frac{7}{2}$, because $c\left(\frac{1}{2} - \frac{1}{p}\right)\geq \frac{35}{4}$.
Thus, in either case,
\[\gamma_{*} \left(\frac{\epsilon}{\lambda}\right)^{-1} \ll_C 2^{-7/2\ln \ln (N)}.\]

Inserting this into (\ref{mark}), we find that
\[\int_{\tilde{B}_{s}}\left| \left| \{a_n \psi_n(x)\}_{n \in I_{k^*,s}} \right|\right|_{V^p}^2 dx \ll_C \ln^{7/2}(N) 2^{-c\ln \ln (N)} 2^{-7/2 \ln \ln (N)} \ll_C 2^{-c \ln \ln (N)}.\]
Now to bound (\ref{V2lastterm}), we apply this to each of the $\leq 2^{c \ln \ln (N)}$ terms, yielding $O(1)$, completing the proof.

\section{Acknowledgements} We thank Mark Rothlisberger for help with translation of related literature.

\texttt{A. Lewko, Department of Computer Science, The University of Texas at Austin}

\textit{alewko@cs.utexas.edu}
\vspace*{0.5cm}

\texttt{M. Lewko, Department of Mathematics, The University of Texas at Austin}

\textit{mlewko@math.utexas.edu}

\end{document}